\documentclass[12pt,reqno]{amsart}
\usepackage{amssymb,amsfonts,amsthm,amsmath,enumitem}
\usepackage{color}
\usepackage{hyperref}
\usepackage[left=1 in,top=1 in,right=1 in,bottom=1 in]{geometry}
\usepackage{cite}

\newcommand{\beq}{\begin{equation}}
\newcommand{\eeq}{\end{equation}}
\newcommand{\beqs}{\begin{equation*}}
\newcommand{\eeqs}{\end{equation*}}
\newcommand{\ba}{\begin{array}}
\newcommand{\ea}{\end{array}}
\newcommand{\beas}{\begin{eqnarray*}}
\newcommand{\eeas}{\end{eqnarray*}}
\newcommand{\bea}{\begin{eqnarray}}
\newcommand{\eea}{\end{eqnarray}}
\newcommand{\bal}{\begin{align}}
\newcommand{\eal}{\end{align}}

\newcommand{\bals}{\begin{align*}}
\newcommand{\eals}{\end{align*}}

\newcommand{\tnum}{\rm(\roman*)}
\newcommand{\rnum}{\rm(\alph*)}

\newcommand{\R}{\ensuremath{\mathbb R}}

\newcommand{\norm}[1]{\| {#1} \|}

\newcommand{\bds}{\begin{displaystyle}}
\newcommand{\eds}{\end{displaystyle}}

\def\longequals{\mathbin{=\kern-2pt=}}
\def\eqdef{\stackrel{\rm def}{=}}

\def\varep{\varepsilon}

\def\ddt{\frac{\rm d}{\rm dt}}
\def\d{{\rm d}}

\newcommand{\remove}[1]{} 
\renewcommand{\remove}[1]{#1} 

\newtheorem{theorem}{Theorem}[section]

\newtheorem{lemma}[theorem]{Lemma}

\newtheorem{proposition}[theorem]{Proposition}

\newtheorem{remark}[theorem]{\bf{Remark}}
\newtheorem{assumption}[theorem]{Assumption}
\theoremstyle{remark}
\newtheorem{example}[theorem]{\bf{Example}}

\newcommand{\essup}{\mathop{\mathrm{ess\,sup}}}

\def\midx{\mu}

\def\myclearpage{}

\definecolor{darkred}{rgb}{.70,.12,.20}

\definecolor{darkgreen}{rgb}{.20,.52,.14}

\numberwithin{equation}{section}

\title{A priori estimates for gaseous flows of Forchheimer-type in heterogeneous porous media}

\author{Emine Celik$^{1}$}
\address{$^{1}$Department of Mathematics, Sakarya University,
54050, Sakarya, T\"{u}rkiye}
\email{eminecelik@sakarya.edu.tr}

\author{Luan Hoang$^{2,*}$}
\address{$^2$Department of Mathematics and Statistics,
Texas Tech University,
1108 Memorial Circle, Lubbock, TX 79409--1042, U. S. A.}
\email{luan.hoang@ttu.edu}

\author{Thinh Kieu$^{3}$}
\address{$^{3}$Department of Mathematics, University of North Georgia, Gainesville Campus,
3820 Mundy Mill Rd., Oakwood, GA 30566, U. S. A.}
\email{thinh.kieu@ung.edu}

\thanks{$^*$Corresponding author.}
\date{\today}



\begin{document}

\begin{abstract} 
We study isentropic fluid flows of gases of the Forchheimer-type in heterogeneous porous media. The governing equation is a  doubly nonlinear parabolic equation with coefficients depending on the spatial variables. Its solutions are subject to a nonlinear Robin boundary condition. We establish the estimates of the solutions for short time in terms of the initial and boundary data. For the proof, the multi-weight versions of the Sobolev inequality, parabolic Sobolev inequality and trace theorem are derived. They are then used to implement the Moser iteration for suitable weighted norms.
\end{abstract}

\maketitle 

\tableofcontents 

\pagestyle{myheadings}\markboth{\sc E. Celik, L. Hoang, and T. Kieu}
{\sc Gaseous flows of Forchheimer-type in heterogeneous porous media
}

\myclearpage    
\section{Introduction and the model}\label{Intro}

We study fluids in porous media with spatial variable $x$, time $t$, velocity $v$, pressure $p$, density $\rho$, porosity $\phi_0\in(0,1)$, absolute viscosity $\mu$ and permeability $k$.
The usual equation of motion for the fluid flow is the Darcy law  \cite{Darcybook}
\beq\label{Darcy}
\frac {\mu}{k} v=-\nabla p + \rho \vec g,
\eeq
where $\vec g$ is the gravitational acceleration.
However, when the Reynolds number is large, there is a deviation from the Darcy law \eqref{Darcy}.
By replacing the linear term in $v$ with nonlinear ones, Forchheimer \cite{Forchh1901,ForchheimerBook} obtained the following three nonlinear models:  
\beq\label{2term}
av+b|v|v=-\nabla p + \rho \vec g,
\eeq
\beq\label{3term}
av+b|v|v+c |v|^2 v=-\nabla p + \rho \vec g,
\eeq
\beq\label{power}
av+d|v|^{m-1}v=-\nabla p + \rho \vec g\text{ for some number } m\in(1, 2).
\eeq
Above, the positive constants $a,b,c,d$ come from experiments.
For more models and discussions about hydrodynamics in porous media, see \cite{MuskatBook,Ward64,BearBook,NieldBook,Scheidegger1974}.

All three Forchheimer  equations \eqref{2term}, \eqref{3term}, \eqref{power} above can be written in a general form  
\beq\label{gF}
\sum_{i=0}^N a_i |v|^{\alpha_i}v=-\nabla p + \rho \vec g.
\eeq  
This is called the generalized Forchheimer equation.
For compressible fluids, by the dimension analysis in \cite{MuskatBook}, equation \eqref{gF} can be modified to become
\beq\label{FM}
\sum_{i=0}^N a_i \rho^{\alpha_i} |v|^{\alpha_i} v=-\nabla p + \rho \vec g.
 \eeq
In particular, Ward \cite{Ward64} established from experimental data the following equation, corresponding to Forchheimer's two term law \eqref{2term},
\beqs
\frac{\mu}{k} v+c_F\frac{\rho}{\sqrt k}|v|v=-\nabla p + \rho \vec g\text{ for some positive constant }c_F.
\eeqs

From the mathematical point of view, the Darcy flows, under the umbrella ``porous medium equations", have been studied intensively  for a long time with vast literature, see e.g. \cite{VazquezPorousBook,AronsonLec1986}. 
However, there is much less mathematical analysis of  the Forchheimer flows. 
For incompressible fluids, see \cite{ChadamQin,Payne1999a,Payne1999b, StraughanBook} and references there in. Regarding compressible fluids, see \cite{ABHI1,HI1,HI2,HIKS1,HKP1,CHIK1,HK1,HK2} for single-phase flows,  \cite{HIK1,HIK2} for two-phase flows, and also \cite{Doug1993,Park2005, Kieu1} for numerical analysis.
In particular, the papers \cite{CH1,CH2} deal with slightly compressible fluids in heterogeneous porous media, while the paper \cite{HK3} deals with anisotropic porous media.

In this article we combine and further develop the techniques in our previous work \cite{CH1,CH2} and \cite{CHK1,CHK2,CHK3,CHK4} to study isentropic fluid flows in heterogeneous porous media.
More specifically, we  consider an open subset $U$ of $\R^n$, for $n=3$ momentarily, which represents the porous media. When the media is heterogeneous, the porosity $\phi_0$ and coefficients $a_i$ in equation \eqref{FM} depend on the spatial variable $x$. Thus, we have $\phi_0=\phi_0(x)$ and  
\beq\label{Fhetero}
 \sum_{i=0}^N a_i(x) \rho^{\bar \alpha_i} |v|^{\bar\alpha_i} v  =-\nabla p + \rho \vec g,
\eeq
where $N\ge 1$, $\bar\alpha_0=0<\bar\alpha_1<\ldots<\bar\alpha_N$ are real numbers, 
each $a_i$ is a function from $\bar U$ to $\R^+\eqdef [0,\infty)$, for $1\le i\le N$, and particularly $a_0(x),a_N(x)>0$ on $\bar U$.

Define the function $g:\bar U\times \mathbb{R}^+\to\R^+$  by
\beq\label{eq2}
g(x,s)=a_0(x)s^{\bar\alpha_0} + a_1(x)s^{\bar\alpha_1}+\cdots +a_N(x)s^{\bar\alpha_N}\quad\text{for } x\in \bar U,\ s\ge 0.
\eeq 
Multiplying equation \eqref{Fhetero} by $\rho$ and using $g(x,s)$ in \eqref{eq2}, we have
\beq\label{eq0}
g(x,\rho|v|)\rho v  =-(\rho \nabla p - \rho^2 \vec g).
\eeq
Denote $Y=\rho \nabla p - \rho^2 \vec g\in\R^n$. Taking the norm of both sides of \eqref{eq0} gives
\beq\label{rvs}
g(x,\rho|v|)\rho|v|=|Y| \text{ which implies } \rho |v|=s(x,|Y|),
\eeq
where $s=s(x,\xi)$, for $x\in \bar U$,  $\xi\in \R^+$ is the unique non-negative solution of  $s g(x,s)=\xi.$
Substituting the last identity  of \eqref{rvs} into $g(x,\rho|v|)$ in \eqref{eq0}, and solving  for $\rho v$ from there, we obtain
\beq\label{rvY}
\rho v=-\frac{Y}{g(x,s(x,|Y|))}.
\eeq
Define the function $K:  \bar U\times \mathbb{R}^+\to \R^+$ by
\beq\label{eq4} K(x,\xi)=\frac{1}{g(x,s(x,\xi))}\text{ for }x\in \bar U,\ \xi\in \R^+,
\eeq
and the vector-valued function $\mathbf X:\bar U\times \R^n\to \R^n$ by 
\beq \label{Xdef}
\mathbf X(x,y)=K(x,|y|)y\text{ for }x\in \bar U,\ y\in\R^n.
\eeq 
The equation \eqref{rvY} reads as 
\beq\label{solve} 
 \rho v=-\mathbf X(x,\rho\nabla p -\rho^2  \vec g).
\eeq
Next, we combine \eqref{solve} with the continuity equation
\beqs
\phi_0\rho_t+{\rm div }(\rho v)=0
\eeqs
to derive
\beq\label{genr}
\phi_0\rho_t={\rm div }\mathbf X(x,\rho\nabla p -\rho^2 \vec g).
\eeq

\medskip\noindent
\emph{Isentropic flows of gases.} The equation of state, in this case, is
$p=c\rho^\gamma$, where $\gamma \ge 1$ is the specific heat ratio.
Setting the pseudo-pressure 
\beq \label{ulam}
u=\frac{c\gamma}{\gamma+1} \rho^{\gamma+1}=\frac{\gamma c^{-1/\gamma}}{\gamma+1}p^\frac{\gamma+1}{\gamma} 
\text{ and letting }\lambda=\frac1{\gamma+1},
\eeq 
we obtain from \eqref{genr} that
\beq \label{ueq}
\phi_0(x)\left(\frac{\gamma+1}{c\gamma}\right)^{\lambda}(u^\lambda)_t
= \nabla \cdot \mathbf X \left (x, \nabla u - \left(\frac{\gamma+1}{c\gamma}u\right)^{2\lambda} \vec g\right).
\eeq

\medskip\noindent
\emph{The boundary condition.} 
Assume that $U$ has a $C^1$-boundary. Let $\vec \nu$ denote the  outward normal vector  on the boundary $\partial U$ of $U$. We  consider the volumetric flux condition $v\cdot \vec\nu=\psi_0(x,t)$ on $\partial U$. This gives 
$ \rho v\cdot \vec\nu =\psi_0 \rho$,
and, thanks to \eqref{solve}, we have
\beq\label{bc1}
\mathbf X \left (x, \nabla u- \left(\frac{\gamma+1}{c\gamma}u\right)^{2\lambda}\vec g\right)\cdot \vec\nu  + \psi_0(x,t) \left(\frac{\gamma+1}{c\gamma}u\right)^\lambda=0 \text{ for }x\in \partial U,\ t>0.
\eeq

\medskip
We  generalize the above equations slightly now. Hereafter, the spatial dimension $n$ is greater or equal to $2$, and $U$ is an open, bounded subset of $\R^n$ with $C^1$-boundary $\Gamma=\partial U$.
Instead of $\lambda$ in \eqref{ulam}, we consider any number $\lambda>0$.
We replace $\left(\frac{\gamma+1}{c\gamma}\right )^{\lambda}\phi_0(x)$ in \eqref{ueq} and $\left(\frac{\gamma+1}{c\gamma}\right)^\lambda\psi_0(x,t)$ in \eqref{bc1} 
with general functions
\beq \label{phps}
\phi:U\to (0,\infty)\text{ and }\psi:\partial U\times(0,\infty)\to\R, \text{ respectively.}
\eeq 
Moreover, the term $-(\frac{\gamma+1}{c\gamma})^{2\lambda}u^{2\lambda} \vec g$ in \eqref{ueq} and \eqref{bc1} is replaced  with $\mathbf Z(u)$, where $\mathbf Z:\R^+\to\R^n$ is a the function that satisfies
\beq\label{Zb}
|\mathbf Z(u)|\le C_Z u^{2\lambda} \text{  for some constant $C_Z>0$ and all $u\ge 0$.}
\eeq

Gathering \eqref{ueq}, \eqref{bc1}, \eqref{phps} and \eqref{Zb}, we study the following initial boundary value problem 
\beq\label{mainpb}
\begin{cases}
\phi(x)(u^\lambda)_t
= \nabla \cdot \mathbf X (x, \nabla u+\mathbf Z(u))&\text{ in }  U\times (0,\infty),\\
u (x,0)=u_0(x) &\text{ on } U,\\
\mathbf X (x,\nabla u+ \mathbf Z(u))\cdot \vec\nu  + \psi(x,t) u^\lambda=0 & \text{ on } \Gamma \times(0,\infty),
\end{cases}
\eeq 
where $u_0(x)$ is a given initial data.
The initial boundary value problem \eqref{mainpb} is the focus of our investigation in this paper. We will derive \emph{a priori} estimates, at least for short time, for its solutions in terms of the initial and boundary data.

The paper is organized as follows.
In Section \ref{Prelim}, we recall the basic properties of the function $\mathbf X(x,y)$ which contains important weight functions $W_1(x)$ and $W_2(x)$, see \eqref{Xy} and \eqref{Xsquare}. Then we establish many inequalities with \emph{multiple} weights, namely, the Sobolev inequality in Lemma \ref{WS1}, the trace theorem in Lemma \ref{trace}, and the parabolic Sobolev inequality in Lemma \ref{WS2}.
It is worth mentioning that the conditions on the weights for the Sobolev inequality 
and trace theorem depend heavily on the exponent $\alpha$. However, for the parabolic Sobolev inequality in Lemma \ref{WS2}, the conditions for the weights are \emph{independent} of $\alpha$. This surprising fact allows us to apply the Moser iteration without imposing increasingly stricter conditions on the weights after each iteration. 
In Section \ref{lasec}, we establish the $L^\alpha$-estimates for a non-negative solution of problem \eqref{mainpb} with sufficiently large $\alpha$, see Theorem \ref{Labound}. The main differential inequality is established in Proposition \ref{Diff4u}.   
In Section \ref{maxsec}, by iterating the inequality for the solution established in Lemma \ref{GLk}, we obtain the main relation between the weighted Lebesgue norm and the $L_{x,t}^\infty$-norm for the solution in Theorem \ref{LinfU}. Then the $L_{x,t}^\infty$-estimate of the solution, for short time, in terms of the initial and boundary data is derived in Theorem \ref{thm45}. 
Thanks to the fixed integral conditions on the weights in the parabolic Sobolev inequality, we can still obtain the $L_{x,t}^\infty$-estimates for the solutions without imposing any undesirable $L^\infty$-conditions for the weights and/or their reciprocals.
As a demonstration, we apply our results to some specific cases  of gas flows in Examples \ref{ex1} and \ref{ex3}.
Finally, the methods developed in this work are applicable to other similar degenerate/singular parabolic equations with the coefficients having rather general  dependence on the spatial variables.
\myclearpage    
\section{Inequalities with multiple weights}\label{Prelim}

First, we recall elementary inequalities that will be used frequently. 
If $x,y\ge 0$, then
\beq\label{ee1}
(x+y)^p\le 2^p(x^p+y^p)\quad  \text{for all }  p>0,
\eeq
\beq\label{ee2}
(x+y)^p\le x^p+y^p\quad  \text{for all } 0<p\le 1,
\eeq
\beq\label{ee3}
(x+y)^p\le 2^{p-1}(x^p+y^p)\quad  \text{for all }  p\ge 1,
\eeq
As a consequence of \eqref{ee3} and the triaangle inequality, one has
\beq\label{ee7}
|x-y|^p\ge 2^{-p+1}|x|^p-|y|^p\quad  \text{for all }  x,y\in\mathbb R^n\text{ and } p\ge 1.
\eeq
We also have the following consequences of Young's inequality. If $x>0$, then 
\beq\label{ee4}
x^\beta \le x^\alpha+x^\gamma\quad \text{for any numbers }\alpha\le \beta\le\gamma,
\eeq
\beq\label{ee5}
x^\beta \le 1+x^\gamma \quad \text{for any numbers } \gamma\ge \beta\ge 0.
\eeq
Clearly, when $x=0$, inequality \eqref{ee4} holds true for $\alpha>0$ and inequality \eqref{ee5} holds true for 
$\beta>0$.

For a Lesbegue measurable function $\varphi:U\to (0,\infty)$, a number $p\ge 1$ and functions $f:U\to \R$, $g:U\times(T_1,T_2)\to \R$ with $T_1<T_2$, we define
\beqs
\|f\|_{L^p_\varphi(U)}=\left(\int_U |f(x)|^p\varphi(x)\d x\right)^{1/p},\quad
\|g\|_{L^p_\varphi(U\times(T_1,T_2))}=\left(\int_{T_1}^{T_2}\int_U |g(x,t)|^p\varphi(x)\d x\d t\right)^{1/p}.
\eeqs

\subsection{Sobolev inequality}\label{Sobsec}
For a  number $p\in[1,n)$, denote by $p^*$ the Sobolev conjugate exponent of $p$, that is, $p^*=np/(n-p)$. 
We recall a convenient interpolation inequality from \cite[Lemma 2.1, Ineq. (2.14)]{CHK1}.
Namely, if $f\in W^{1,p}(U)$ with $p\in[1,n)$ then
\beq\label{scare}
\|f\|_{L^{p^*}(U)}\le c_1\|\nabla f\|_{L^p(U)}+c_2\|f\|_{L^1(U)}.
\eeq
where $c_1$ and $c_2$ are positive constants depending on $U$, $n$ and $p$.
Since \cite{CHK1} did not contain a proof for \eqref{scare}, we quickly verify it here.

\begin{proof}[Proof of Ineq. \eqref{scare}]
We have the standard Sobolev inequality
\beq \label{basic}
\|f\|_{L^{p^*}(U)}\le c(\|\nabla f\|_{L^p(U)}+\|f\|_{L^p(U)})\text{ for some constant } c>0.
\eeq 

If $p=1$, then \eqref{basic} implies \eqref{scare}.
Consider $p>1$. Then $1<p<p^*$ and there is $\theta\in(0,1)$ such that $1/p=\theta/1 +(1-\theta)/p^*$.
By the interpolation inequality for Lebesgue integrals and Yong's inequality, we have 
\beq\label{cf}
c\|f\|_{L^p(U)}\le c\|f\|_{L^1(U)}^\theta \|f\|_{L^{p^*}(U)}^{1-\theta} \le \frac12 \|f\|_{L^{p^*}(U)} + 2^{(1-\theta)/\theta} c^{1/\theta}\|f\|_{L^1(U)}.
\eeq
Using \eqref{cf} to estimate the last term in \eqref{basic} and absorbing the $L^{p^*}$-norm of $f$ to the left-hand side, we obtain \eqref{scare}.
\end{proof}

We now establish a suitable Sobolev inequality with multiple weights.

\begin{lemma}\label{WS1}
Let numbers $p>1$ and $r_1>0$ satisfy 
\beq\label{rone}
\frac{n}{n+p}<r_1<1\le r_1p<n,
\eeq 
and denote
\beq\label{rstar}
r_*=\frac{nr_1+r_1p-n}{nr_1}= 1+\frac p n-\frac1{r_1}\in(0,1).
\eeq
Given numbers $r\ge 0$ and $s \in\R$. Assume $\alpha\in\R$ satisfies
\begin{align} \label{manyalpha}
&\alpha\ge s,\quad \alpha > \frac{p-s}{p-1},\\
\label{many12}
& \alpha> \frac{2(r+s-p)}{r_*}.
\end{align} 
Denote
\beq\label{powers_def}
\begin{cases}
\begin{aligned}
m&= \frac{\alpha-s+p}{p}\in[1,\alpha), \\
\theta&=\frac{\alpha+2r}{\alpha(1+r_*)+2(p-s)}\in(0,1),\\
 \mu_1&=\frac{r+\theta(s-p)}{1-\theta}\in (-\alpha,\infty).
 \end{aligned}
\end{cases}
\eeq
Let $\omega(x)\ge 0$, $\varphi(x)>0$, $W(x)>0$ be functions defined on $U$.
Define
\beq \label{Gdef}
\begin{aligned}
    G_1&=\left(\int_U \varphi(x)^{-\frac {\alpha-s+p}{\alpha(p-1)+s-p}} \d x\right)^\frac{\alpha(p-1)+s-p}{\alpha},\quad 
        G_2=\left(\int_U W(x)^{-\frac{r_1}{1-r_1}}\d x\right)^{\frac{1-r_1}{r_1}},  \\
  G_3&=\left(\int_U  \varphi(x)^{-1}\omega(x)^{\frac{1}{(1-\theta)(1+\mu_1/\alpha)}}\d x\right)^{1+\mu_1/\alpha},\quad
  \Phi_1= G_1^\theta G_3^{1-\theta},\quad 
   \Phi_2=  G_2^{\frac{ \theta}{1-\theta}} G_3.
\end{aligned}
\eeq 
Then  one has, for any function $u\in W^{1,{r_1pm}}(U)$ and number $\varepsilon>0$, that
\beq \label{S10}
\begin{aligned}
&\int_U |u(x)|^{\alpha+r} \omega(x) \d x 
 \le \varep \int_U |u(x)|^{\alpha-s}|\nabla u(x)|^p W(x)\d x \\
 &\quad +D_{1,m,\theta} \Phi_1 \left(\int_U |u(x)|^\alpha\varphi(x)\d x \right)^{1+r/\alpha}
 +\varep^{-\frac\theta{1-\theta}}D_{2,m,\theta} \Phi_2 \left( \int_U |u(x)|^\alpha \varphi(x)\d x \right)^{1+\mu_1/\alpha},
\end{aligned}
\eeq 
where 
\beq \label{dthe}
D_{1,z,\eta}=(c_4 2^{z})^{\eta p} ,\quad 
D_{2,z,\eta}=(c_3z 2^{z})^\frac{\eta p}{1-\eta}
\text{ for $z>0$, $\eta\in(0,1)$}
\eeq 
with $c_3$, respectively, $c_4$ being the positive number $c_1$, respectively, $c_2$ in \eqref{scare} for $p:=r_1p$.
\end{lemma}
\begin{proof} 
Note that from the first and last inequalities in \eqref{rone}, one has $r_*>0$ and $r_*<1$, respectively, thus verifying the last property of $r_*$ in \eqref{rstar}.
Regarding the number $\alpha$, by considering $s>1$ and $s\le 1$ in \eqref{manyalpha}, we have
\beq\label{ag1}
\alpha>1.
\eeq
It also follows from \eqref{manyalpha} that $1\le m<\alpha$, which verifies the stated property of $m$ in \eqref{powers_def}. 
Set
 \beq\label{bardef1}
 \overline{p}=r_1 p,\quad \overline{\alpha}=r_1\alpha, \quad\overline{s}=r_1 s,\quad 
 \overline{m}=\frac{\overline{\alpha}-\overline{s}+\overline{p}}{\overline{p}} .
 \eeq
By \eqref{rone}, we have $1\le \overline{p}<n$, hence can set  $\overline{q}=\overline {p}^*\overline{ m}.$
Obviously,
 \beq\label{bardef2}
 \overline{m}=\frac{\alpha-s+p}{p}=m\in[1,\alpha),\text{ and } 
 \overline{q}=(r_1 p)^* m=\frac{nr_1(\alpha-s+p)}{n-r_1p}.
 \eeq

Because $r\ge 0$ and $\alpha > 0$, we have $\alpha/2<\alpha+r$.
Thanks to the fact $r_1>\frac{n}{n+p}$ and $r_1p<n$ in \eqref{rone}, we deduce that $r_1p+nr_1-n>0$ and $n-r_1p>0$.
Then imposing the condition
\beq\label{alphaq}
\alpha+r<\overline{q}
\eeq
is equivalent to requiring
\beq\label{alplarge0}
\begin{aligned}
\alpha
&>\frac{nr_1(s-p)+r(n-r_1p)}{r_1p+nr_1-n}
=\frac1{r_*}\left\{ s-p+\frac{r}{r_1n}(n-r_1p-nr_1+nr_1)\right\}\\
&=\frac{s-p+r(1-r_*)}{r_*}=\frac{r+s-p}{r_*}-r.     
\end{aligned}
\eeq
If $r+s-p\le 0$, then the last condition \eqref{alplarge0} is trivially met thanks to our consideration  $\alpha>0$ and $r\ge 0$.
If $r+s-p>0$, then it is met thanks to  \eqref{many12}.
In both cases, the requirement \eqref{alplarge0} is met, and hence we have \eqref{alphaq}.

  Now that  $\alpha/2<\alpha+r<\overline q$, we set
  \beq \label{inval}
  \text{$\widetilde\theta$ and $\theta_0$ to be two numbers in the interval $(0,1)$}
  \eeq such that
\beq \label{thethe}
\alpha+r=(1-\widetilde \theta) \alpha/2 +\widetilde\theta \overline q\text{ and } \frac1{\alpha+r}=\frac{1-\theta_0}{\alpha/2}+\frac{\theta_0}{\overline q}.
\eeq 
Explicitly, we have
\beq \label{prethe}
\widetilde\theta=\frac{(\alpha+r)-(\alpha/2)}{\overline q-\alpha/2}
\text{ and } 
\theta_0=\frac{1/(\alpha+r)-1/(\alpha/2)}{1/\overline q-1/(\alpha/2)}
=\frac{1/(\alpha/2)-1/(\alpha+r)}{1/(\alpha/2)-1/\overline q}.
\eeq 
On the one hand, we can calculate  $\theta_0$ and use the formula of $\overline q$ in \eqref{bardef2} to obtain
\beqs
\theta_0= \frac{\overline{q}(\alpha+2r)}{(\alpha+r)(2\overline{q}-\alpha)}
=\frac{nr_1(\alpha-s+p)(\alpha+2r)}{(\alpha+r)[\alpha(2nr_1+r_1p-n)+2nr_1(p-s)]},    
\eeqs
which yields
\beq\label{thezero}
\theta_0=\frac{(\alpha-s+p)(\alpha+2r)}{(\alpha+r)[\alpha(1+r_*) +2(p-s)]}.
\eeq 
On the other hand, we calculate the numerator and denominator of $\theta_0$ in \eqref{prethe} separately  and rewrite it in the form
$$\theta_0
=\left(\frac{(\alpha+r)-\alpha/2}{(\alpha+r)\cdot \alpha/2}\right)
\left(\frac{\overline q-\alpha/2}{\bar q \cdot \alpha/2} \right)^{-1}
=\frac{(\alpha+r)-\alpha/2}{\overline q-\alpha/2} \cdot \frac{\overline q}{\alpha+r}=\widetilde\theta \cdot \frac{\overline q}{\alpha+r},$$
hence,
\beq\label{tthe1}
\widetilde\theta \overline q=\theta_0 (\alpha+r).
\eeq 
By switching the roles $\alpha/2\leftrightarrow \overline q$, $\widetilde\theta \leftrightarrow (1-\widetilde\theta)$ and  $\theta_0 \leftrightarrow (1-\theta_0)$,  we have  from \eqref{tthe1} that
\beq \label{tthe2}
(1-\widetilde\theta)(\alpha/2)=(1-\theta_0)(\alpha+r).
\eeq 

Using the first equation in \eqref{thethe} to rewrite $|u|^{\alpha+r}$ and then applying  H\"older's inequality, we have
\begin{align*}
 \int_U |u|^{\alpha+r}\omega \d x 
 &=\int_U |u|^{\widetilde \theta \overline q } \cdot\left (|u|^{(1-\widetilde \theta)\alpha/2}\omega\right) \d x 
 \le \left(\int_U |u|^{\overline q}\d x \right)^{\widetilde\theta}
 \left(\int_U |u|^{\alpha/2} \omega^{1/(1-\widetilde\theta)}\d x \right)^{1-\widetilde\theta}.
\end{align*}
The last integral can be estimated by the Cauchy--Schwarz inequality as follows
\begin{align*}
 &\int_U |u|^{\alpha/2} \omega^{1/(1-\widetilde\theta)}\d x 
 = \int_U \left(|u|^{\alpha/2}\varphi^{1/2}\right)\cdot\left(\varphi^{-1/2} \omega^{1/(1-\widetilde\theta)} \right)\d x \\
 & \le \left(\int_U |u|^\alpha \varphi \d x  \right)^{1/2}
 \left(\int_U  \varphi^{-1}\omega^{2/(1-\widetilde\theta)}\d x \right)^{1/2}
=\|u\|_{L^\alpha_\varphi(U)}^{\alpha/2} \left(\int_U  \varphi^{-1}\omega^{2/(1-\widetilde\theta)}\d x \right)^{1/2}.
\end{align*}
Thus,
\begin{align*}
 \int_U |u|^{\alpha+r}\omega \d x 
 & \le \| u\|_{L^{\overline{q}}(U)}^{\widetilde \theta \overline q}
 \|u\|_{L^\alpha_\varphi(U)}^\frac{(1-\widetilde \theta)\alpha}{2}
 \left(\int_U  \varphi^{-1}\omega^{2/(1-\widetilde\theta)}\d x \right)^\frac{1-\widetilde\theta}{2}.
\end{align*}

By using formulas of $\widetilde \theta$ in \eqref{tthe1} and of $(1-\widetilde\theta)$ in \eqref{tthe2}, we rewrite the preceding estimate as
\beq\label{maine}
 \int_U |u|^{\alpha+r}\omega \d x 
  \le R_1\| u\|_{L^{\overline{q}}(U)}^{\theta_0 (\alpha+r)}\|u\|_{L^\alpha_\varphi(U)}^{(1-\theta_0)(\alpha+r)} ,
\eeq
where
\beq\label{R1}
R_1=\left(\int_U  \varphi^{-1}\omega^{2/(1-\widetilde\theta)}\d x \right)^\frac{1-\widetilde \theta}{2}
=\left(\int_U  \varphi^{-1}\omega^{\frac{\alpha}{(1-\theta_0) (\alpha+r)}}\d x \right)^\frac{(1-\theta_0)(\alpha+r)}{\alpha}.
\eeq 

We estimate  the norm $\| u\|_{L^{\overline{q}}}$ in \eqref{maine}.
Since $\overline m\ge 1$ and $1\le\overline{p}<n$, applying inequality \eqref{scare} to function $f=|u|^{\overline m}$ and power $p:=\overline p$  yields
\beq\label{sobwm}
\begin{aligned}
\|\,|u|^{\overline{m}}\,\|_{L^{\overline{p}^*}(U)} \le c_3 \|\nabla (|u|^{\overline{m}})\|_{L^{\overline{p}}(U)} +c_4\|\,|u|^{\overline{m}}\,\|_{L^1(U)}.
\end{aligned}
\eeq
Note from the definition \eqref{bardef1} of $\overline m$ that $(\overline{m}-1)\overline{p}=\overline{\alpha}-\overline{s}$. Hence \eqref{sobwm} can be written as
\beqs
\Big(\int_U |u|^{\overline{p}^*\overline{m}} \d x \Big)^{1/\overline{p}^*}
\le c_3 \overline{m}\Big(\int_U |u|^{\overline{\alpha}-\overline{s}}|\nabla u|^{\overline{p}} \d x \Big)^{1/\overline{p}} +c_4 \int_U |u|^{\overline{m}} \d x .
\eeqs
Raising both sides to the power $1/\overline{m}\le 1$ and using inequality \eqref{ee2},  we obtain
\beq\label{t1}
\Big(\int_U |u|^{\overline{p}^*\overline{m}} \d x \Big)^{1/(\overline{p}^*\overline{m})}
\le (c_3 \overline{m})^{1/\overline{m}}\Big(\int_U |u|^{\overline{\alpha}-\overline{s}}|\nabla u|^{\overline{p}} \d x \Big)^{1/(\overline{p}\, \overline{m})} +c_4^{1/\overline{m}} \Big(\int_U |u|^{\overline{m}} \d x \Big)^{1/\overline{m}}.
\eeq
Then \eqref{t1} yields
\beq\label{So1}
\|u\|_{L^{\overline{q}}(U)}\le (c_3 m)^\frac 1m \Big(\int_U |u|^{\overline{\alpha}-\overline{s}}|\nabla u|^{\overline{p}} \d x \Big)^\frac1{\overline{\alpha}-\overline{s}+\overline{p}} +c_4^\frac1m\|u\|_{L^m(U)}.
\eeq
For the integral with the gradient term in \eqref{So1}, we rewrite 
$$|u|^{\overline{\alpha}-\overline{s}}|\nabla u|^{\overline{p}} =\left ( |u|^{\overline{\alpha}-\overline{s}}|\nabla u|^{\overline{p}} W^{r_1}\right )\cdot W^{-r_1},$$
and apply H\"older's inequality with powers $1/r_1>1$ and $1/(1-r_1)$ noticing that $\overline{\alpha}-\overline{s}+\overline{p}=r_1(\alpha-s+p)$. 
It results in
 \begin{equation}\label{229a}
  \|u\|_{L^{\overline{q}}(U)}\le (c_3m)^\frac{1}{m} \Big(\int_U |u|^{\alpha-s}|\nabla u|^p W\d x \Big)^\frac{1}{\alpha-s+p}\cdot R_2+c_4^\frac 1{m}\|u\|_{L^ m(U)},
 \end{equation}
where
\beq\label{R2}
R_2=\Big(\int_U W(x)^{-\frac{r_1}{1-r_1}}\d x \Big)^{\frac{1-r_1}{r_1(\alpha-s+p)}}.
\eeq 
Now using \eqref{229a} in \eqref{maine} together with inequality \eqref{ee1} applied to the power  $p:=\theta_0(\alpha+r)$, we obtain
\beq\label{ualpr}
\begin{aligned}
 \int_U |u|^{\alpha+r}\omega \d x 
 \le 2^{\theta_0 (\alpha+r)}  R_1 &\Big\{  (c_3m)^\frac{\theta_0 (\alpha+r)}{m} \Big(\int_U |u|^{\alpha-s}|\nabla u|^p W\d x \Big)^\frac{\theta_0(\alpha+r)}{\alpha-s+p}  R_2^{\theta_0 (\alpha+r)}\\
 &\quad +c_4^\frac{\theta_0 (\alpha+r)}{m}\|u\|_{L^m(U)}^{\theta_0 (\alpha+r)}\Big\} 
  \|u\|_{L_\varphi^{\alpha}(U)}^{(1-\theta_0)(\alpha+r)}.
\end{aligned}
\eeq
Denote
\beq\label{thetadef}
\theta = \frac{\theta_0(\alpha+r)}{\alpha-s+p}.
\eeq
With the use of formula \eqref{thezero} for $\theta_0$, this number $\theta$  actually assumes the same value given in \eqref{powers_def}.
Then \eqref{ualpr} can be written as
\beq\label{ualpr2}
\begin{aligned}
 &\int_U |u|^{\alpha+r}\omega \d x \\
 &\le \left(\varep\int_U |u|^{\alpha-s}|\nabla u|^p W\d x \right)^\theta \cdot \left(\varep^{-\theta} 2^{\theta_0 (\alpha+r)}(c_3m)^\frac{\theta_0 (\alpha+r)}{m}  R_1  R_2^{\theta_0(\alpha+r)}
 \|u\|_{L_\varphi^{\alpha}(U)}^{(1-\theta_0)(\alpha+r)} \right)\\
 &\quad + 2^{\theta_0 (\alpha+r)}c_4^\frac{\theta_0 (\alpha+r)}{m} R_1 \|u\|_{L^m(U)}^{\theta_0 (\alpha+r)} \|u\|_{L_\varphi^{\alpha}(U)}^{(1-\theta_0)(\alpha+r)}.
\end{aligned}
\eeq
It is clear from \eqref{thetadef} and the fact $\alpha\ge s$ that $\theta>0$.
Now, with $\theta$ found by the formula in \eqref{powers_def}, the condition $\theta<1$ is equivalent to 
$\displaystyle \alpha > 2(r-p+s)/r_*$ which was already assumed in \eqref{many12}. 
Thus, $\theta$ belongs to the interval $(0,1)$ as stated in \eqref{powers_def}.
Applying Young's inequality to the first term on the right-hand side of \eqref{ualpr2} with powers   $1/\theta$ and $1/(1-\theta)$, we obtain
\beq\label{S2}
\begin{aligned}
 \int_U |u|^{\alpha+r}\omega(x) \d x 
 & \le  \varep \int_U |u|^{\alpha-s}|\nabla u|^p W(x)\d x \\
&\quad + \varep^{-\frac{\theta}{1-\theta}} 2^\frac{\theta_0 (\alpha+r)}{1-\theta} (c_3{m})^\frac{\theta_0(\alpha+r)}{{m}(1-\theta)} R_1^\frac1{1-\theta}R_2^\frac{\theta_0(\alpha+r)}{1-\theta}  \|u\|_{L_\varphi^{\alpha}(U)}^\frac{(1-\theta_0)(\alpha+r)}{1-\theta}\\
&\quad 
+ 2^{\theta_0 (\alpha+r)}c_4^\frac{\theta_0 (\alpha+r)}{m} \cdot R_1  \|u\|_{L^m(U)}^{\theta_0 (\alpha+r)} \|u\|_{L_\varphi^{\alpha}(U)}^{(1-\theta_0)(\alpha+r)}.
\end{aligned}
\eeq
Since $m< \alpha$, applying H\"older's inequality with the power $\alpha/m$ and $\alpha/(\alpha-m)$ to bound the $L^m$-norm of $u$ on the right-hand side of \eqref{S2}, we obtain
\begin{align*}
    \|u\|_{L^m(U)}^{\theta_0(\alpha+r)}=\left( \int_U \left( |u|^m \varphi^{m/\alpha} \right) \cdot \varphi^{-m/\alpha} \d x \right)^{\theta_0(\alpha+r)/m}
    &\le \|u\|_{L_\varphi^\alpha(U)}^{\theta_0(\alpha+r)}\cdot R_3,
\end{align*}
where 
\beq\label{R3}
R_3= \Big(\int_U \varphi^{-\frac m{\alpha-m}} \d x \Big)^{(\frac 1 m -\frac 1 \alpha)\theta_0(\alpha+r)}.
\eeq
Thus, we obtain
\begin{multline}\label{S11}
\int_U |u|^{\alpha+r}\omega(x) \d x 
\le \varep \int_U |u|^{\alpha-s}|\nabla u|^p W(x)\d x \\
+ \varep^{-\frac{\theta}{1-\theta}} 2^\frac{\theta_0 (\alpha+r)}{1-\theta} (c_3{m})^\frac{\theta_0(\alpha+r)}{{m}(1-\theta)}R_1^\frac1{1-\theta} R_2^\frac{\theta_0(\alpha+r)}{1-\theta}  \|u\|_{L_\varphi^{\alpha}(U)}^\frac{(1-\theta_0)(\alpha+r)}{1-\theta} 
 +2^{\theta_0 (\alpha+r)}c_4^\frac{\theta_0 (\alpha+r)}{m} R_1 R_3\|u\|_{L_\varphi^\alpha(U)}^{\alpha+r}.
\end{multline}

We recalculate many powers in inequality \eqref{S11} to express them in terms of $\theta$ instead of $\theta_0$.
From \eqref{thetadef}, we have 
\begin{align}\label{thep}
\theta_0(\alpha+r)&=\theta(\alpha-s+p), \quad \frac{\theta_0(\alpha+r)}{m}=\theta p,\\
\frac{(1-\theta_0)(\alpha+r)}{1-\theta}
&=\frac{(\alpha+r)-\theta(\alpha-s+p)}{1-\theta}=\frac{(1-\theta)\alpha+r+\theta(s-p)}{1-\theta}
=\alpha+\midx_1.\label{mub}
\end{align}
Thanks to \eqref{mub} and the fact $\theta_0,\theta\in(0,1)$, we must have $\alpha+\mu_1>0$, which implies $\mu_1>-\alpha$ as stated in \eqref{powers_def}.
We use the last expression in \eqref{mub} to conveniently rewrite
$$
\frac{(1-\theta_0)(\alpha+r)}{\alpha}=(1-\theta)(1+\midx_1/\alpha).
$$
Therefore, $R_1$ in \eqref{R1} becomes
$R_1=G_3^{1-\theta}$.
For $R_2$ in \eqref{R2}, it is clear that 
$R_2=G_2^\frac{1}{\alpha-s+p}$,  hence, with the first relation in \eqref{thep}, one has 
$$R_2^\frac{\theta_0(\alpha+r)}{1-\theta}
=G_2^{\frac{1}{\alpha-s+p}\cdot \frac{\theta(\alpha-s+p)}{1-\theta}}=G_2^\frac{\theta}{1-\theta}.
$$
For $R_3$ in \eqref{R3}, we calculate
\begin{align*}
\frac{m}{\alpha-m}&=\frac{\alpha-s+p}{\alpha(p-1)+s-p},\\
\theta_0 (\alpha+r)\Big(\frac 1m- \frac 1\alpha\Big)
&=\theta(\alpha-s+p)\Big(\frac{p}{\alpha-s+p}-\frac1\alpha\Big)
=\frac{\theta(\alpha(p-1)+s-p)}{\alpha}.
\end{align*}
Thus,
$R_3=G_1^\theta.$
We can rewrite the terms involving $R_1$, $R_2$, $R_3$ in \eqref{S11} as
\beqs
R_1^\frac1{1-\theta} R_2^\frac{\theta_0(\alpha+r)}{1-\theta} =G_3G_2^\frac{\theta}{1-\theta}=\Phi_2
\text{ and } 
R_1 R_3=G_3^{1-\theta} G_1^\theta=\Phi_1.
\eeqs

Regarding the constants in \eqref{S11}, by using the second relation in \eqref{thep}, we have
$$  2^\frac{\theta_0 (\alpha+r)}{1-\theta} (c_3{m})^\frac{\theta_0(\alpha+r)}{{m}(1-\theta)}
= (2^mc_3{m})^\frac{\theta_0(\alpha+r)}{{m}(1-\theta)}=(2^mc_3{m})^\frac{\theta p}{1-\theta}=D_{2,m,\theta},$$
$$ 2^{\theta_0 (\alpha+r)}c_4^\frac{\theta_0 (\alpha+r)}{m} 
= (2^m c_4)^\frac{\theta_0 (\alpha+r)}{m}=(2^m c_4)^{\theta p}=D_{1,m,\theta}.$$

Combining these calculations with \eqref{S11}, we obtain 
\begin{align*}
\int_U |u|^{\alpha+r}\omega(x) \d x 
&\le \varep \int_U |u|^{\alpha-s}|\nabla u|^p W(x)\d x \\
&\quad + \varep^{-\frac{\theta}{1-\theta}} D_{2,m,\theta}\Phi_2  \|u\|_{L_\varphi^{\alpha}}^{\alpha+\mu_1}
 +D_{1,m,\theta} \Phi_1 \|u\|_{L_\varphi^\alpha}^{\alpha+r}
\end{align*}
which implies the desired inequality  \eqref{S10}.
\end{proof}

\subsection{Trace theorem} \label{tracesec}
We first establish a simple trace theorem with one weight function which will assist us in the proof of the main trace theorem -- Lemma \ref{trace} below.

\begin{lemma}\label{wtrace}
If  $p>1$, $s\in\mathbb{R}$ and 
\beq\label{alphac1}
\alpha \ge \max \left\{s, \frac{p-s}{p-1}\right\},
\eeq 
then one has, for any function $u\in W^{1,\alpha}(U)$ and number $\varepsilon>0$, that
\beq\label{trace10}
\begin{split}
    \int_\Gamma |u(x)|^\alpha \d S 
&\le \varepsilon \int_U |u|^{\alpha-s}|\nabla u(x)|^p W(x)\d x  + c_5 \int_U |u(x)|^\alpha \d x   \\
&\quad  + (c_6 \alpha)^\frac p{p-1} \varepsilon^{-\frac 1{p-1}} \int_U |u(x)|^{\alpha+\frac{s-p}{p-1}} W(x)^{-\frac 1{p-1}}\d x ,
\end{split}
\eeq
where   $c_5,c_6>0$ are constants depending on $U$, but not on $u(x),\alpha,s,p$, see \eqref{firstrace} below. 
\end{lemma}
\begin{proof}
By considering $s\ge 1$ and $s<1$ in \eqref{alphac1}, we have $\alpha\ge 1$.
We recall the trace theorem
\beq\label{firstrace}
\int_\Gamma |f|\d S \le c_5\int_U |f|\d x +c_6\int_U |\nabla f|\d x ,
\eeq
for all $f \in W^{1,1}(U)$, where $c_5$ and $c_6$ are positive constants depending on $U$.
Applying this trace theorem to $f=|u|^{\alpha}$, we have 
\beqs
\int_\Gamma |u|^\alpha \d S \le c_5 \int_U |u|^\alpha \d x  + c_6\alpha \int_U |u|^{\alpha-1}|\nabla u|\d x.
\eeqs
Rewriting $c_6\alpha |u|^{\alpha-1}|\nabla u|$ in the last integral as a product of 
$$(\varepsilon W)^{1/p}|u|^{\frac{\alpha-s}p}|\nabla u|\text{ and }c_6\alpha (\varepsilon W)^{-1/p} |u|^\frac{(p-1)\alpha+s-p}{p},$$
and applying Young's inequality with powers $p$ and $p/(p-1)$,  we obtain  inequality \eqref{trace10}.
\end{proof}

Our main trace theorem with two weights is the following.

\begin{lemma}\label{trace}
Let numbers  $p>1$ and  $r_1>0$ satisfy \eqref{rone}, and the number $r_*$ be defined by \eqref{rstar}.
Given numbers $s\in \R$ and  $r\ge 0$.
 Assume $\alpha$ is a number that satisfies \eqref{manyalpha}.
Let the numbers $m$, $\theta$ and $\midx_1$ be defined as in \eqref{powers_def}, and set
\beq\label{rtilde} 
\widetilde r=r+\frac{r+s-p}{p-1}=\frac{rp+s-p}{p-1}.
\eeq
Let $c_5, c_6$ be as in Lemma \ref{wtrace} and 
 $D_{1,z,\eta}, D_{2,z,\eta}$ as in \eqref{dthe}.
Let $\varphi(x)>0$ and $W(x)>0$ be functions defined on $U$.
Define  $G_1$, $G_2$ as in \eqref{Gdef}, and 
\beq\label{Ph34}
    G_4=\left(\int_U  \varphi(x)^{-1} \d x\right)^{1+\mu_1/\alpha},\quad 
    \Phi_3=G_1^\theta G_4^{1-\theta}, \quad
    \Phi_4=G_2^{\frac{ \theta}{1-\theta}} G_4.
\eeq

Let $\varep>0$ and  $u(x)$ be any function on $U$ satisfy $|u|^\alpha\in W^{1,1}(U)$ and $|u|^m\in W^{1,p}(U)$.

\begin{enumerate}[label=\tnum]
    \item\label{trneg} Case $\widetilde r <0$. Suppose $\alpha>|\widetilde r|$.
Let 
\beq \label{Ph5}   
  \Phi_5=\left(\int_U W(x)^{\frac {\alpha}{(p-1)\widetilde r}}\varphi(x)^{\frac{\alpha+\widetilde r}{\widetilde r}}(x)\d x\right)^{-\widetilde r/\alpha}.
\eeq
Then one has
\beq \label{trace0}
\begin{aligned}
&\int_\Gamma |u(x)|^{\alpha+r} \d S 
\le 2\varepsilon \int_U |u(x)|^{\alpha-s}|\nabla u(x)|^p W(x) \d x
 +z_1 \Phi_3 \left(\int_U |u(x)|^\alpha\varphi(x) \d x \right)^{1+ r/\alpha}\\
&\quad +\varep^{-\frac\theta{1-\theta}} z_2 \Phi_4\left(\int_U |u(x)|^\alpha\varphi(x) \d x \right)^{1+\mu_1/\alpha}
  + \varepsilon^{-\frac 1{p-1}}  z_3 \Phi_5 \left(\int_U |u(x)|^\alpha\varphi(x) \d x \right)^{1+ \widetilde r/\alpha},
\end{aligned}
\eeq 
where 
\beq \label{zzz}
z_1=c_5 D_{1,m,\theta},\quad 
z_2= c_5^{\frac{1}{1-\theta}}D_{2,m,\theta},\quad 
z_3= (c_6 (\alpha+r))^\frac p{p-1}.
\eeq 

\item\label{trpos} Case $\widetilde r \ge 0$.
Suppose 
\beq\label{only}  
\alpha> \frac{2(\widetilde r+s-p)}{r_*} =\frac{2p(r+s-p)}{r_*(p-1)}.
\eeq 
Let 
\beq\label{tmtilde} 
\widetilde\theta=\frac{\alpha+2\widetilde{r}}{\alpha(1+r_*)+2(p-s)}\in(0,1), \quad 
 \widetilde\mu_1=\frac{\widetilde r+\widetilde\theta(s-p)}{1-\widetilde\theta}\in(-\alpha,\infty),
\eeq
\beq \label{Ph67}
G_5=\left(\int_U  \varphi(x)^{-1}W(x)^{-\frac{1}{(p-1)(1-\widetilde\theta) (1+\widetilde\mu_1/\alpha)}}\d x\right)^{1+\widetilde\mu_1/\alpha},\  
\Phi_6=  G_1^{\widetilde\theta} G_5^{1-\widetilde\theta}, \ 
\Phi_7=  G_2^{\frac{\widetilde\theta}{1-\widetilde\theta}} G_5.
\eeq
Then one has
\begin{align} \label{trace1}
&\int_\Gamma |u(x)|^{\alpha+r} \d S 
\le 3\varepsilon \int_U |u(x)|^{\alpha-s}|\nabla u(x)|^p W(x) \d x 
+ z_1 \Phi_3 \left(\int_U |u(x)|^\alpha\varphi(x) \d x \right)^{1+r/\alpha} \notag \\
&\quad + \varep^{-\frac\theta{1-\theta}}z_2 \Phi_4 \left(\int_U |u(x)|^\alpha\varphi(x) \d x \right)^{1+\mu_1/\alpha}
+ \varepsilon^{-\frac 1{p-1}}z_4 \Phi_6
  \left( \int_U |u(x)|^\alpha \varphi(x) \d x \right)^{1+\widetilde r/\alpha}\notag \\
&\quad + \varep^{-(\frac 1 {p-1}+\frac p{p-1}\cdot \frac{\widetilde\theta}{1-\widetilde\theta} )}z_5 \Phi_7
  \left( \int_U |u(x)|^\alpha \varphi(x) \d x \right)^{1+\widetilde\mu_1/\alpha},
\end{align}
where $z_1$, $z_2$, $z_3$ are from \eqref{zzz}, $z_4=z_3 D_{1,m,\widetilde\theta}$ and $z_5=z_3^\frac1{1-\widetilde\theta} D_{2,m,\widetilde\theta}$.
\end{enumerate}
\end{lemma}
\begin{proof}
Since $\alpha$ satisfies \eqref{manyalpha}, we recall from \eqref{ag1} that $\alpha>1$.
Below, we denote
\beqs
I_0=\int_U  W|u|^{\alpha-s}|\nabla u|^p \d x
\text{ and } I_1= \int_U |u|^\alpha\varphi \d x.
\eeqs 
Applying inequality \eqref{trace10} in Lemma \ref{wtrace} to the powers $\alpha:=\alpha+r$ and $s:=s+r$ yields
\beq\label{rep}
\begin{aligned}
\int_\Gamma |u|^{\alpha+r} \d S 
&\le \varepsilon I_0 + c_5 \int_U |u|^{\alpha+r} \d x   + C_0 \varepsilon^{-\frac 1{p-1}} \int_U W^{-\frac 1{p-1}}|u|^{\alpha+\widetilde r} \d x,
\end{aligned}
\eeq
where $C_0=(c_6 (\alpha+r))^\frac p{p-1}$ which is $z_3$ in \eqref{zzz}.
The conditions in \eqref{alphac1} for $\alpha:=\alpha+r$ and $s:=s+r$  become
\beqs
\alpha\ge s,\quad \alpha\ge \frac{p-s-rp}{p-1},
\eeqs
which, in fact, are met thanks to our assumption \eqref{manyalpha}.

In order to estimate the second term on the right-hand side of \eqref{rep}, we apply  inequality \eqref{S10}  in Lemma \ref{WS1} with $\omega(x)=1$ and parameter $\varep:=\varep c_5^{-1} $. Note that $G_3$, $\Phi_1$, $\Phi_2$ in \eqref{Gdef} now become $G_4$, $\Phi_3$, $\Phi_4$ in \eqref{Ph34}, respectively. It results in 
\beq \label{term1}
c_5\int_U |u|^{\alpha+r} \d x
  \le \varep I_0 
  +c_5 D_{1,m,\theta}\Phi_3   I_1^{1+r/\alpha}
  +\varep^{-\frac\theta{1-\theta}}c_5^{\frac 1{1-\theta}}D_{2,m,\theta} \Phi_4 I_1^{1+\mu_1/\alpha}.
\eeq 

Inequality \eqref{term1} requires condition \eqref{manyalpha}, which we already assumed, and condition \eqref{many12}, which we rewrite here
\begin{align}\label{tal1}
  \alpha> \frac{2(r+s-p)}{r_*}.
\end{align}
If $r+s-p\le 0$, then condition  \eqref{tal1} is met thanks to the simple fact $\alpha>0$.
If $r+s-p>0$, then \eqref{tal1}  follows from assumption \eqref{only} and the fact $p/(p-1)>1$.

For the third term on the right-hand side of \eqref{rep}, we consider the following two cases of $\widetilde{r}$.

\medskip\noindent\textit{Case $\widetilde{r}<0.$} 
In this case, we must have $r+s-p<0$, hence condition \eqref{tal1} is met.
 Assume $\alpha>|\widetilde r|$. Since $\alpha+\widetilde{r}>0$, we can estimate the  integral $\displaystyle\int_U W^{-\frac 1{p-1}}|u|^{\alpha+\widetilde r}\d x$ in \eqref{rep}
by  H\"{o}lder's inequality with powers $-\alpha/\widetilde r$ and  $\alpha/(\alpha+\widetilde r)$ as follows, with the use of $\Phi_5$ in \eqref{Ph5},
\beq\label{term2a}
\begin{aligned}
&    \int_U W^{-\frac 1{p-1}}|u|^{\alpha+\widetilde r} \d x
=\int_U  (W^{-\frac 1{p-1}}\varphi^{-\frac{\alpha+\widetilde r}{\alpha}})\cdot (\varphi^{\frac{\alpha+\widetilde r}{\alpha}}|u|^{\alpha+\widetilde r})\d x\\
&\le \left(\int_U W^{\frac {\alpha}{(p-1)\widetilde r}}\varphi^{\frac{\alpha+\widetilde r}{\widetilde r}}\d x \right)^{-\frac{\widetilde r}{\alpha}} \left(\int_U |u|^\alpha\varphi \d x\right)^{\frac{\alpha+\widetilde r}{\alpha}} =\Phi_5 I_1^{\frac{\alpha+\widetilde r}{\alpha}}.
\end{aligned}
\eeq
Combining \eqref{rep},  \eqref{term1} and \eqref{term2a},  we have
\beqs
\int_\Gamma |u|^{\alpha+r} \d S 
\le 2\varepsilon I_0 
  +c_5 D_{1,m,\theta}\Phi_3   I_1^{1+r/\alpha}
  +\varep^{-\frac\theta{1-\theta}}c_5^{\frac 1{1-\theta}}D_{2,m,\theta} \Phi_4 I_1^{1+\mu_1/\alpha}
  + C_0 \varepsilon^{-\frac 1{p-1}} \Phi_5I_1^{1+\widetilde r/\alpha}.
\eeqs
This proves \eqref{trace0}.

\medskip\noindent\textit{Case $\widetilde{r}\ge 0$.} 
Given any number $\varepsilon_1>0$. Estimating the last integral in \eqref{rep}
now by using \eqref{S10} in Lemma \ref{WS1} with the weight $\omega(x) = W(x)^{-\frac 1{p-1}}$ this time and replacing $r$ by $\widetilde{r}$, we obtain
\beq \label{lastint}
\begin{aligned}
\int_U W^{-\frac 1{p-1}}|u|^{\alpha+\widetilde r} \d x
  &\le \varep_1 I_0
  +D_{1,m,\widetilde \theta} \Phi_6   I_1^{1+\widetilde r/\alpha}
   +\varep_1^{-\frac{\widetilde\theta}{1-\widetilde\theta}}D_{2,m,\widetilde \theta} \Phi_7 I_1^{1+\widetilde\mu_1/\alpha},
\end{aligned}
\eeq 
where, from  \eqref{powers_def} with $r:=\widetilde r$, the powers $\widetilde \theta$ and $\widetilde \mu_1$ are given by \eqref{tmtilde},
and $\Phi_6$ and $\Phi_7$ are defined as in \eqref{Ph67}.
The first requirement \eqref{manyalpha} is independent of $r$ and was already assumed. The second requirement \eqref{many12} for $r:=\widetilde r$, noting that $r_1,p,s,r_*$ are not changed,  becomes 
\beqs
 \alpha> \frac{2(\widetilde{r}+s-p)}{r_*}=\frac{2p(r+s-p)}{r_*(p-1)}.
\eeqs 
This condition, in fact, is already assumed in \eqref{only}. 

Because $\widetilde\theta$ and $\widetilde\mu_1$ in \eqref{lastint} play the role of $\theta$ and $\mu$ in \eqref{S10}, we have from the properties of $\theta$ and $\mu$ in \eqref{powers_def} that $\widetilde\theta\in(0,1)$ and $\widetilde\mu_1\in(-\alpha,\infty)$ as stated in \eqref{tmtilde}.

Multiplying \eqref{lastint} by $B:=C_0 \varepsilon^{-\frac 1{p-1}} $ and choosing 
$\varep_1:=B^{-1}\varep=\varep^{\frac p{p-1}}C_0^{-1}$
result in 
\beq \label{term2}
\begin{aligned}
&C_0\varepsilon^{-\frac 1{p-1}}\int_U W^{-\frac 1{p-1}}|u|^{\alpha+\widetilde r} \d x\\
&\le \varepsilon I_0 
  +C_0\varepsilon^{-\frac 1{p-1}}D_{1,m,\widetilde \theta}  \Phi_6   I_1^{1+\widetilde r/\alpha}
+C_0\varepsilon^{-\frac 1{p-1}}  (\varep^{\frac p{p-1}}C_0^{-1})^{-\frac{\widetilde\theta}{1-\widetilde\theta}}  D_{2,m,\widetilde \theta}\Phi_7 I_1^{1+\widetilde\mu_1/\alpha}
\\
&  =\varepsilon I_0 
 + C_0\varepsilon^{-\frac 1{p-1}}D_{1,m,\widetilde \theta} \Phi_6   I_1^{1+\widetilde r/\alpha}
+ \varep^{-(\frac 1 {p-1}+\frac p{p-1}\frac{\widetilde\theta}{1-\widetilde\theta} )}C_0^\frac 1{1-\widetilde\theta} D_{2,m,\widetilde \theta} \Phi_7 I_1^{1+\widetilde\mu_1/\alpha}
 .
\end{aligned}
\eeq 
Combining estimates \eqref{term1} and \eqref{term2} with inequality \eqref{rep}, we derive
\begin{align*}
\int_\Gamma |u|^{\alpha+r} \d S
&\le 3\varepsilon I_0
+c_5 D_{1,m,\theta} \Phi_3   I_1^{1+r/\alpha}
+ \varep^{-\frac\theta{1-\theta}}c_5^{\frac 1{1-\theta}} D_{2,m,\theta}\Phi_4 I_1^{1+\mu_1/\alpha} \\
&\quad 
 + C_0\varepsilon^{-\frac 1{p-1}}D_{1,m,\widetilde \theta} \Phi_6   I_1^{1+\widetilde r/\alpha}
+ \varep^{-(\frac 1 {p-1}+\frac p{p-1}\frac{\widetilde\theta}{1-\widetilde\theta} )}C_0^\frac 1{1-\widetilde\theta} D_{2,m,\widetilde \theta} \Phi_7 I_1^{1+\widetilde\mu_1/\alpha}
 ,
\end{align*}
which proves \eqref{trace1}.
\end{proof}

 \subsection{Parabolic Sobolev inequality}\label{paraSsec}
We derive a parabolic Sobolev inequality that contains multiple weights.

\begin{lemma}\label{WS2}
Let numbers  $p>1$ and  $r_1>0$ satisfy \eqref{rone}, and the number $r_*$ be defined by \eqref{rstar}.
Let $\varphi(x)>0$ and $W(x)>0$ be functions defined on $U$.
Define
\beq \label{Jss}
 \Phi_8=\left ( \int_U  \varphi(x)^{\frac{2}{r_*}-1}\d x\right )^\frac{r_*}{2}\left[ \left(\int_U W(x)^{-\frac{r_1}{1-r_1}}\d x \right)^{1-r_1}
+\left(\int_U \varphi^{-\frac{r_1}{1-r_1}}\d x\right)^{1-r_1}\right]^\frac{1}{r_1}.
\eeq 
Given a number $s\in \R$, assume $\alpha$ is a number such that 
\beq \label{alcond}
    \alpha\ge s,\ 
    \alpha> \frac{p-s}{p-1} \text{ and }
    \alpha >\frac{2(s-p)}{r_*}.
\eeq
Let 
\begin{align}\label{defkappa}
 \kappa& = 1+ \frac{r_*}{2}+\frac{p-s}{\alpha}\in (1,\infty), \\ 
\label{tz}
\theta_0 &=\frac 1{1+\frac{r_*\alpha}{2(\alpha-s+p)} }\in (0,1),\\
\label{mdef}
m&=\frac{\alpha-s+p}{p}\in[1,\alpha).
\end{align}
Let $u(x,t)$ be any function on $U\times (0,T)$ such that for almost every $t\in (0,T)$ the function $x\in U\mapsto u(x,t)$ belongs to $W^{1,r_1pm}(U)$.
Then one has
\beq\label{ppsi1}
\begin{aligned}
\|u\|_{L^{\kappa \alpha}_\varphi(U\times(0,T))}
&\le  (c_7^p m^\frac1{r_1} \Phi_8)^\frac{1}{\kappa\alpha}  \Bigg(\int_0^T\int_U |u(x,t)|^{\alpha-s}|\nabla u(x,t)|^p W(x)\d x\d t \\
&\quad +\int_0^T\int_U |u(x,t)|^{\alpha-s+p} \varphi(x)  \d x\d t\Bigg)^\frac{1}{\kappa\alpha}
 \cdot  \essup_{t\in(0,T)} \|u(\cdot,t)\|_{L_\varphi^{\alpha}(U)}^{1-\theta_0},
\end{aligned}
\eeq
where the positive constant $c_7$ is independent of $\alpha$, $T$ and $u(x,t)$, see \eqref{stdSov} below.
\end{lemma}
\begin{proof}
With the first two conditions in \eqref{alcond}, we have \eqref{manyalpha}. 
Thus, recalling \eqref{ag1} and \eqref{powers_def}, we further have $\alpha>1$ and $m\in[1,\alpha)$.
We neglect formulas \eqref{defkappa} and \eqref{tz} for $\kappa$ and $\theta_0$ momentarily.
Let $\kappa$  be any number greater than $1$. 
Let $\overline p$, $\overline m$ and $\overline q$ be as in \eqref{bardef1} and \eqref{bardef2}.
We apply inequality \eqref{maine} in the proof of Lemma \ref{WS1} with 
\beq \label{rtemp}
r=(\kappa-1)\alpha  \text{ and } \omega=\varphi
\eeq 
so that $\alpha+r=\kappa \alpha$.
The conditions for \eqref{maine} are 
\beq \label{pra2}
p>1,\ \eqref{rone},\ r\ge 0,\  \eqref{manyalpha}\text{ and }\eqref{alplarge0}.
\eeq 
Note in \eqref{pra2} that the first, second and fourth conditions are already met, while the third condition comes from the formula of $r$ in \eqref{rtemp} and the facts $\alpha>0$ and $\kappa>1$. The last  requirement \eqref{alplarge0} reads as
\beq\label{ar1}
\begin{aligned}
\alpha >\frac{s-p+r(1-r_*)}{r_*}=\frac{s-p+(\kappa-1)\alpha(1-r_*)}{r_*}  .  
\end{aligned}
\eeq
Also, the number $\theta_0$ from \eqref{inval} and \eqref{thezero} takes the particular form
\beq\label{th1}
\begin{aligned}
\theta_0&=\frac{(\alpha-s+p)(\alpha+2r)}{(\alpha+r)[\alpha(1+r_*) +2(p-s)]}=\frac{(\alpha-s+p)(2\kappa-1)}{\kappa[\alpha(1+r_*) +2(p-s)]}.    
\end{aligned}
\eeq 
Then applying inequality \eqref{maine} to the function $u=u(\cdot,t)$ for almost every $t\in(0,T)$, we have 
\beq\label{maine2}
 \int_U |u(x,t)|^{\kappa\alpha}\varphi(x) \d x 
  \le K_1\| u(\cdot,t)\|_{L^{\overline{q}}(U)}^{\theta_0 \kappa\alpha}\|u(\cdot,t)\|_{L^\alpha_\varphi(U)}^{(1-\theta_0)\kappa\alpha} ,
\eeq
where $K_1$ plays the role of $R_1$ in \eqref{maine} and is defined by, see \eqref{R1} with $\omega=\varphi$, 
\beq\label{K1} 
K_1=\left(\int_U  \varphi^{-1}\varphi^{\frac{\alpha}{(1-\theta_0) (\kappa\alpha)}}\d x\right)^\frac{(1-\theta_0)(\kappa\alpha)}{\alpha}=\left(\int_U  \varphi^{\frac{1}{(1-\theta_0) \kappa}-1}\d x\right)^{(1-\theta_0)\kappa}.
\eeq 

In the following calculations up to \eqref{preineq2}, we  temporarily denote $u=u(\cdot,t)$.
We estimate  the norm $\| u\|_{L^{\overline{q}}(U)}$ in \eqref{maine2}.
Because $1\le\overline{p}<n$, we have the standard Sobolev's inequality 
\beq \label{stdSov}
\|f\|_{L^{\overline p^*}(U)} \le c_7 \left( \int_U |\nabla f|^{\overline p}\d x + \int_U |f|^{\overline p}\d x  \right)^{1/\overline p}
\text{ for any function $f\in W^{1,\overline p}(U)$,}
\eeq
where $c_7$ is a positive constant.
With $\overline m\ge 1$, applying inequality \eqref{stdSov} to function $f=|u|^{\overline m}$  yields
\beq\label{stdum}
\begin{aligned}
\|\,|u|^{\overline{m}}\,\|_{L^{\overline{p}^*}(U)}  
\le c_7 \left( \overline m \int_U |u|^{(\overline m-1)\overline p}|\nabla u|^{\overline p}\d x + \int_U |u|^{\overline  m \,\overline p}\d x  \right)^{1/\overline p}.
\end{aligned}
\eeq
Note from the definition \eqref{bardef1} of $\overline m$ that $(\overline{m}-1)\overline{p}=\overline{\alpha}-\overline{s}$. Hence \eqref{stdum} can be written as
\beqs
\left(\int_U |u|^{\overline{p}^*\overline{m}} \d x \right)^{1/\overline{p}^*}
\le c_7 \left (\overline{m}\int_U |u|^{\overline{\alpha}-\overline{s}}|\nabla u|^{\overline{p}} \d x + \int_U |u|^{\overline  m \,\overline p} \d x \right)^{1/\overline{p}}.
\eeqs
Raising both sides to the power $1/\overline{m}$ and using the facts
$ \overline{m}\, \overline{p}^*=\overline q$,
$ \overline{m}\, \overline{p}=\overline{\alpha}-\overline{s}+\overline{p}$, $\overline{m}=m$, we obtain
\beq\label{stduq}
\|u\|_{L^{\overline{q}}(U)}
=\|u\|_{L^{\overline{m}\, \overline{p}^*}(U)}
\le 
 c_7^{1/m} \left (m\int_U |u|^{\overline{\alpha}-\overline{s}}|\nabla u|^{\overline{p}} \d x + \int_U |u|^{\overline{\alpha}-\overline{s}+\overline{p}} \d x \right)^\frac1{\overline{\alpha}-\overline{s}+\overline{p}}.
\eeq
For the last two integrals in \eqref{stduq}, we rewrite 
$$|u|^{\overline{\alpha}-\overline{s}}|\nabla u|^{\overline{p}} =\left ( |u|^{\overline{\alpha}-\overline{s}}|\nabla u|^{\overline{p}} W^{r_1}\right )\cdot W^{-r_1},\quad 
|u|^{\overline{\alpha}-\overline{s}+\overline{p}}=(|u|^{\overline{\alpha}-\overline{s}+\overline{p}}\varphi^{r_1})\cdot \varphi^{-r_1}.
$$
We apply H\"older's inequality with powers $1/r_1>1$ and $1/(1-r_1)$ to each integral noticing that $\overline{\alpha}-\overline{s}+\overline{p}=r_1(\alpha-s+p)$. 
It results in
 \beqs
  \|u\|_{L^{\overline{q}}(U)}\le c_7^\frac{1}{m}\left[ m \left(\int_U |u|^{\alpha-s}|\nabla u|^p W\d x \right)^{r_1}\cdot K_2
  +\left(\int_U |u|^{\alpha-s+p} \varphi \d x\right)^{r_1}\cdot K_3 \right]^\frac{1}{r_1(\alpha-s+p)},
 \eeqs
where
\beq\label{K23}
K_2=\Big(\int_U W(x)^{-\frac{r_1}{1-r_1}}\d x \Big)^{1-r_1},\quad 
K_3=\left(\int_U \varphi^{-\frac{r_1}{1-r_1}}\d x\right)^{1-r_1}.
\eeq 
Consequently,
\begin{equation}\label{2b}
\begin{aligned}
  \|u\|_{L^{\overline{q}}(U)}
  &\le c_7^\frac{1}{m}\left[  
  \left(\int_U |u|^{\alpha-s}|\nabla u|^p W\d x +\int_U |u|^{\alpha-s+p}  \varphi \d x\right)^{r_1} (m K_2+ K_3) \right]^\frac{1}{r_1(\alpha-s+p)}\\
&\le c_7^\frac{1}{m}  [m (K_2+ K_3)]^\frac{1}{r_1(\alpha-s+p)} \left(\int_U |u|^{\alpha-s}|\nabla u|^p W\d x +\int_U |u|^{\alpha-s+p}  \varphi \d x\right)^\frac{1}{\alpha-s+p} .
\end{aligned}
\end{equation}

Now raising \eqref{2b} to the power $\theta_0\kappa\alpha$ and combining it with \eqref{maine2}, we obtain 
\beq\label{preineq}
 \int_U |u|^{\kappa\alpha}\varphi \d x 
 \le  c_7^\frac{\theta_0 \kappa\alpha}{m} m^\frac{\theta_0 \kappa\alpha}{r_1(\alpha-s+p)} K_4  \left(\int_U |u|^{\alpha-s}|\nabla u|^p W\d x +\int_U |u|^{\alpha-s+p}  \varphi \d x\right)^\frac{\theta_0 \kappa\alpha}{\alpha-s+p} 
  \|u\|_{L_\varphi^{\alpha}(U)}^{(1-\theta_0)\kappa\alpha},
\eeq
where 
\beq \label{K4}
K_4=K_1 (K_2+ K_3)^\frac{\theta_0 \kappa\alpha}{r_1(\alpha-s+p)}.
\eeq 
We impose that
\beq\label{teq2}
\frac{\theta_0 \kappa \alpha}{\alpha-s+p}=1.
\eeq
Under  assumption \eqref{teq2},  one has $\theta_0\kappa\alpha/m=(\alpha-s+p)/m=p.$
Thus, inequality \eqref{preineq} becomes
\beq\label{preineq2}
\begin{aligned}
 \int_U |u|^{\kappa\alpha}\varphi \d x
  &\le  c_7^{p}m^\frac1{r_1} K_4  \left(\int_U |u|^{\alpha-s}|\nabla u|^p W\d x +\int_U |u|^{\alpha-s+p} \varphi  \d x\right)
  \|u\|_{L_\varphi^{\alpha}(U)}^{(1-\theta_0)\kappa\alpha}.
\end{aligned}
\eeq
Integrating \eqref{preineq2} in $t$ from $0$ to $T$ and taking the essential supremum in $t$ of the last term give
\beq\label{preineq3}
\begin{aligned}
\int_0^T\int_U |u|^{\kappa \alpha}\varphi(x) \d x \d t
& \le  c_7^{p} m^\frac1{r_1} K_4  \left(\int_0^T\int_U |u|^{\alpha-s}|\nabla u|^p W\d x\d t +\int_U |u|^{\alpha-s+p} \varphi  \d x\d t\right)\\ &\quad \times \essup_{t\in (0,T)} \|u(\cdot,t)\|_{L_\varphi^{\alpha}(U)}^{(1-\theta_0)\kappa \alpha}.
\end{aligned}
\eeq

Now using \eqref{th1} and \eqref{teq2}, we have 
\beqs
1=\frac{\theta_0 \kappa \alpha}{\alpha-s+p}=\frac{(2\kappa-1)\alpha}{\alpha(1+r_*) +2(p-s)},
\eeqs
which yields 
\beq\label{samekap}
\kappa=\frac{2+r_*}{2}+\frac{p-s}{\alpha}=1+\frac{r_*}{2}+\frac{p-s}{\alpha}.
\eeq
This formula of $\kappa$ is, in fact, the same one stated in \eqref{defkappa}.
With \eqref{samekap}, on the one hand, the condition $\kappa>1$ is equivalent to  
\beq\label{rela}
\alpha>\frac{2(s-p)}{r_*},
\eeq
which we already assumed in \eqref{alcond}.
On the other hand, condition \eqref{ar1} now becomes 
\beqs
\alpha >\frac{s-p+\left(\frac{r_*}{2}+\frac{p-s}\alpha \right)\alpha (1-r_*)}{r_*}=\frac{\alpha(1-r_*)}{2}+s-p,
\eeqs
which is equivalent to  
\beq\label{simal}
\alpha>\frac{2(s-p)}{r_*+1}.
\eeq
Note that $\alpha>0$ implies \eqref{simal}  in the case  $s-p\le 0$, and \eqref{rela} implies \eqref{simal} in the case  $s-p>0$. Therefore, condition \eqref{alcond} already covers the requirement \eqref{simal}. 

Solving for $\theta_0$ from \eqref{teq2} and using formula \eqref{samekap} for $\kappa$ yield 
\begin{align*}
    \theta_0=\frac{\alpha-s+p}{\kappa\alpha}=\frac{\alpha-s+p}{\alpha-s+p+ \alpha r_*/2}=\frac 1{1+\frac{r_*\alpha}{2(\alpha-s+p)} }, 
\end{align*}
which is the formula \eqref{tz} for $\theta_0$. The fact $\theta_0\in(0,1)$ in \eqref{tz} comes from the original requirement \eqref{inval}. 

We calculate the power $(1-\theta_0)\kappa$ in formula \eqref{K1} of $K_1$. By the formula \eqref{defkappa} of $\kappa$ and relation \eqref{teq2}, we have
\beq\label{tk}
(1-\theta_0)\kappa=\kappa-\kappa\theta_0
=\left( 1+\frac{p-s}{\alpha}+\frac{r_*}{2}\right)-\kappa\theta_0=\frac{\alpha+p-s}{\alpha}+\frac{r_*}{2}-\kappa\theta_0=\frac{r_*}{2}.
\eeq
Hence, $K_1= ( \int_U  \varphi(x)^{\frac{2}{r_*}-1}\d x )^\frac{r_*}{2}$.
By this identity, the formulas of $K_2$, $K_3$ in \eqref{K23}, and by using \eqref{teq2} for the power of $(K_2+K_3)$ in the formula \eqref{K4} of $K_4$, we have 
\beqs 
K_4= K_1
 (K_2+ K_3)^\frac{1}{r_1}=\Phi_8.
\eeqs 

Finally, taking the power $(\kappa\alpha)^{-1}$ of \eqref{preineq3}, we obtain \eqref{ppsi1}.
\end{proof}

\section{Integral estimates}\label{lasec}
We study the initial boundary value problem \eqref{mainpb}.
With the power $\bar\alpha_N$ in \eqref{eq2}, we denote and use throughout this paper the number
\beq\label{eq9}
a=\frac{\bar\alpha_N}{\bar\alpha_N+1}\in (0,1).
\eeq

We recall some  properties of $K(x,\xi)$ defined by \eqref{eq4}.
Define the main weight functions
\beqs 
M_*(x)=\max\{ a_j(x):j=0,\ldots, N\},\quad
m_*(x)=\min\{a_0(x),a_N(x)\},
\eeqs
\beqs
W_1(x)=\frac{a_N(x)^a}{2 N M_*(x)},\quad W_2(x)= \frac{NM_*(x)}{a_N(x)^{1-a}m_*(x)}.
\eeqs
From Lemma 1.1 of  \cite{CH1}, we have for all $\xi\ge 0$ that
\beqs 
\frac{2W_1(x)}{\xi^a +a_N(x)^a}\le K(x,\xi)\le \frac{W_2(x)}{\xi^a},
\eeqs
\beqs 
W_1(x)\xi^{2-a}-\frac{a_N(x)}{2} \le K(x,\xi)\xi^2\le W_2(x)\xi^{2-a}.
\eeqs
Consequently, we have the following main properties of the function $\mathbf X(x,y)$ in \eqref{Xdef}, for $x\in \bar U$  and $y\in\R^n$,
\beq \label{Xy}
|\mathbf X(x,y)|=K(x,|y|)|y| \le W_2(x)|y|^{1-a},
\eeq
\beq \label{Xsquare}
W_1(x)|y|^{2-a}-\frac{a_N(x)}{2} \le \mathbf X(x,y)\cdot y=K(x,|y|)|y|^2\le W_2(x)|y|^{2-a}.
\eeq

We will also use the following combined weight
\beq \label{W3}
W_3(x) =W_1(x)+\frac{W_2(x)^{2-a}}{W_1(x)^{1-a}}=W_1(x)+\frac{2^{1-a}NM_*(x)}{a_N(x)^{2(1-a)}m_*(x)^{2-a}}.
\eeq 
In Proposition \ref{Diff4u} below, we will apply Lemmas \ref{WS1} and \ref{trace} to the parameters 
$p=2-a$ and $s=\lambda+1$.
Then condition \eqref{rone} for the number $r_1$ reads as
\beqs
\frac{n}{n+2-a}<r_1<1\le r_1(2-a)<n.
\eeqs
For $r_1<1$ and $n\ge 2$, the last inequality is already satified. Hence, condition \eqref{rone} is reduced to
\beq\label{newro}
\frac{n}{n+2-a}<r_1<1\le r_1(2-a).
\eeq 
Also the number $r_*$ defined by \eqref{rstar} now becomes 
\beq\label{newrs}
r_*= 1+\frac{2-a}n-\frac1{r_1}\in(0,1).
\eeq

We will use the notation $\psi^-=\max\{0,-\psi\}$. Then $\psi^-\ge 0$ and $-\psi\le \psi^-$.

Below $u(x,t)$ is a non-negative solution of the problem \eqref{mainpb}. We assume that it has enough regularity so that all our calculations are valid.
We aim to estimate the $L^\alpha_\phi$-norm of $u(x,t)$ at least for small time $t>0$ .

\begin{proposition}\label{Diff4u}
Given a number $r_1$ satisfying \eqref{newro}, denote the number $r_*$ by \eqref{newrs}.
Let $r$ be a fixed number such that  
\beq \label{newrr}
r>\max\{0, \lambda(3-2a)-1\},
\eeq
and denote
\beq \label{newtr}
 \widetilde r=r+\frac{r-1+\lambda+a}{1-a}.
\eeq 
Consider  a number $\alpha$ satisfying
\beq\label{newaa}
\alpha> \lambda+1 \text{ and } \alpha>\frac{2(2-a)(r+a+\lambda-1)}{r_*(1-a)}.
\eeq
Let
\beq\label{mtm}
\begin{aligned}
\theta&=\Theta(r)\eqdef\frac{\alpha+2r}{\alpha(1+r_*)+2(1-a-\lambda)}\in(0,1), \\ 
 \mu_1&= \Lambda(r,\theta)\eqdef  \frac{r+\theta(a+\lambda-1)}{1-\theta}\in(-\alpha,\infty),
\end{aligned}
\eeq
and 
\beq \label{trtm2}
\widetilde\theta=\Theta(\widetilde r)\in(0,1),\quad 
\widetilde\mu_1=\Lambda(\widetilde r,\widetilde \theta)\in(-\alpha,\infty).
\eeq
Assume that 
\begin{align}
{\mathcal K}_1&\eqdef \int_U  \phi(x)^{-1}\d x,\quad 
{\mathcal K}_2\eqdef \int_U \phi(x)^{-\frac {\alpha+1-\lambda-a}{\alpha(1-a)-1+\lambda+a}} \d x,\quad 
{\mathcal K}_3\eqdef \int_U a_N(x)^{\frac{\alpha}{\lambda+1} }\phi(x)^{1-\frac{\alpha}{\lambda+1} } \d x \notag\\
{\mathcal K}_4&\eqdef \int_U W_1(x)^{-\frac{r_1}{1-r_1}}\d x,\quad 
{\mathcal K}_5\eqdef \int_U W_3(x)^{\frac{\alpha+r}{r+1-\lambda(3-2a)}} \d x, \label{KKs}\\
 {\mathcal K}_6&\eqdef \int_U  \phi(x)^{-1}(x)W_1(x)^{-1/[(1-a)(1-\widetilde\theta) (1+\widetilde\mu_1/\alpha)]}\d x \notag
\end{align}
 are finite numbers.  Set
\beq\label{mmm2}
\mu_{\min}= \max\{\lambda+1,-\mu_1,-\widetilde\mu_1\}\in(0,\alpha),\quad \mu_{\max}=\max\{\mu_1,\widetilde r,\widetilde\mu_1\}>0.
\eeq
Then there is a number $Z_*>0$ independent of $u_0(x)$, $\psi(x,t)$ and $u(x,t)$ such that, for $t>0$,
\beq\label{diffIneq}
\begin{aligned}
&\ddt\int_U \phi(x) u^\alpha(x,t) \d x
+ \frac{\alpha(\alpha-\lambda)}{2^{3-a}\lambda} \int_U u^{\alpha-\lambda-1}(x,t) |\nabla u(x,t)|^{2-a} W_1(x)\d x\\
&\le Z_*\left\{ \left( \int_U u^\alpha(x,t) \phi(x)\d x\right)^{1-\frac{\mu_{\min}}\alpha}+\left( \int_U u^\alpha(x,t) \phi(x)\d x\right)^{1+\frac{\mu_{\max}}\alpha}+ M(t)\right\},
\end{aligned}
\eeq
where 
\beq\label{Mtdef}
M(t) =1+\int_\Gamma (\psi^-(x,t))^{\frac{\alpha+r}{r}} \d S.
\eeq
\end{proposition}
\begin{proof}
Let $c_3$ and $c_4$ be the positive constants $c_1$ and $c_2$, respectively, in \eqref{scare} with $p:=r_1(2-a)$. Define $D_{1,z,\eta}$ and $D_{2,z,\eta}$ by in \eqref{dthe} with $p:=2-a$. For the sake of convenience in later calculations, we list here the correspondence between the numbers $G_j$ in Lemmas \ref{WS1} and \ref{trace} when
\beq\label{choice}
\omega(x)=1,\quad W(x)=W_1(x),\quad \varphi(x)=\phi(x),\quad  p=2-a,\quad s=1+\lambda,
\eeq 
with the numbers ${\mathcal K}_j$ in the statement of this proposition
\beq\label{GKs}
    G_1={\mathcal K}_2^\frac{\alpha(1-a)-1+\lambda+a}{\alpha}  ,\quad G_2={\mathcal K}_4^\frac{1-r_1}{r_1},\quad G_3=G_4={\mathcal K}_1^{1+\mu_1/\alpha},\quad G_5={\mathcal K}_6^{1+\widetilde\mu_1/\alpha}.
\eeq

Define two quantities 
$$ I_0(t)=\int_U |u(x,t)|^{\alpha-\lambda-1}|\nabla u(x,t)|^{2-a} W_1(x)\d x
\text{ and }
J_0(t)=\int_U u(x,t)^\alpha\phi(x) \d x.$$
Multiplying the PDE in \eqref{mainpb} by $u^{\alpha-\lambda}$, integrating over domain $U$, and using integration by parts, we obtain 
\begin{align*}
&\frac{\lambda}{\alpha}\ddt\int_U \phi u^\alpha \d x
=\frac{\lambda}{\alpha}\ddt\int_U \phi (u^\lambda)^{\alpha/\lambda} \d x
=\int_U \phi(x)(u^\lambda)_t u^{\alpha-\lambda} \d x
\\
&= -(\alpha-\lambda)\int_U  \mathbf X (x,\nabla u+\mathbf Z(u))\cdot(\nabla u)  u^{\alpha-\lambda-1} \d x
+\int_\Gamma  \mathbf X (x,\nabla u+\mathbf Z(u))\cdot \vec \nu u^{\alpha-\lambda} \d S.
\end{align*}
Utilizing the boundary condition in \eqref{mainpb} for the last surface integral, we obtain
\beq\label{dJ1}
\frac\lambda\alpha J_0' = -(\alpha-\lambda) I_1 -\int_\Gamma  \psi u^\alpha  \d S \le -(\alpha-\lambda) I_1 +I_2,
\eeq
where
\beqs
I_1=\int_U  \mathbf X (x,\nabla u+\mathbf Z(u))\cdot(\nabla u)  u^{\alpha-\lambda-1} \d x
\text{ and }
I_2=\int_\Gamma  \psi^- u^\alpha  \d S.
\eeqs

Rewrite the second gradient $\nabla u$ in $I_1$ as $(\nabla u+\mathbf Z(u) )-\mathbf Z(u)$, we have 
\beq \label{Ionea}
I_1=\int_U \mathbf X (x,\nabla u+\mathbf Z(u))\cdot(\nabla u+\mathbf Z(u))  u^{\alpha-\lambda-1} \d x 
-\int_U  \mathbf X (x,\nabla u+\mathbf Z(u))\cdot  \mathbf Z(u)\, u^{\alpha-\lambda-1} \d x.
\eeq 

For the first integral in \eqref{Ionea}, we use the first inequality of \eqref{Xsquare} and inequality \eqref{ee7} with $p=2-a$ to have
\begin{align*}
    \mathbf X (x,\nabla u+\mathbf Z(u))\cdot(\nabla u+\mathbf Z(u)) 
    &\ge  W_1|\nabla u+\mathbf Z(u)|^{2-a}  - \frac 1 2 a_N \\
    &\ge  W_1 (2^{a-1}|\nabla u|^{2-a}-|\mathbf Z(u)|^{2-a})  - \frac 1 2 a_N .
\end{align*}
Denote $c_0=2^{a-1}\in(0,1)$. Combining this with the estimate of $\mathbf Z(u)$ in \eqref{Zb}, we obtain
\beq\label{uX1}
\begin{aligned}
&u^{\alpha-\lambda-1}  \mathbf X (x,\nabla u+\mathbf Z(u))\cdot(\nabla u+\mathbf Z(u)) \\
&\ge  u^{\alpha-\lambda-1} W_1\Big(c_0|\nabla u|^{2-a}  -C_Z^{2-a} u^{2\lambda(2-a)} \Big) - \frac 1 2 u^{\alpha-\lambda-1}a_N \\
&=  c_0 u^{\alpha-\lambda-1} W_1|\nabla u|^{2-a} 
- C_Z^{2-a} W_1 u^{\alpha+\lambda(3-2a)-1}  - \frac 1 2 u^{\alpha-\lambda-1}a_N .
\end{aligned}
\eeq 

For the second integral in \eqref{Ionea},  applying inequalities   \eqref{Xy} and \eqref{ee2} gives
\begin{align*}
| \mathbf X (x,\nabla u+\mathbf Z(u))\cdot  \mathbf Z(u)| u^{\alpha-\lambda-1} 
&\le W_2 |\nabla u+\mathbf Z(u)|^{1-a}  |\mathbf Z(u)| u^{\alpha-\lambda-1}\\
& \le W_2 (|\nabla u|^{1-a}+|\mathbf Z(u)|^{1-a})|\mathbf Z(u)| u^{\alpha-\lambda-1}\\
&\le C_Z W_2 |\nabla u|^{1-a} u^{(\alpha-\lambda-1)+2\lambda} + C_Z^{2-a} W_2u^{\alpha+\lambda(3-2a)-1}.
\end{align*}
Given a number $\varep_1>0$, rewriting the second to last term as an appropriate product and applying Young's inequality with powers $(2-a)/(1-a)$ and $(2-a)$ give
\begin{align*}
&C_Z W_2 |\nabla u|^{1-a} u^{(\alpha-\lambda-1)+2\lambda} \\
&=\left\{ (\varep_1  W_1 u^{\alpha-\lambda-1})^{(1-a)/(2-a)}|\nabla u|^{1-a} \right\} \cdot 
\left\{ C_Z(\varep_1  W_1)^{-(1-a)/(2-a)}W_2 u^{(\alpha-\lambda-1)/(2-a)+2\lambda}\right\}\\
&\le \varep_1 W_1 |\nabla u|^{2-a}u^{\alpha-\lambda-1}
 + C_Z^{2-a} \varep_1^{-(1-a)} W_1^{-(1-a)} W_2^{2-a} u^{\alpha+\lambda(3-2a)-1} . 
\end{align*}
Thus,
\beq\label{uX2}
\begin{aligned}
&| \mathbf X (x,\nabla u+\mathbf Z(u))\cdot  \mathbf Z(u)| u^{\alpha-\lambda-1} 
\le \varep_1 W_1 |\nabla u|^{2-a}u^{\alpha-\lambda-1}\\
&\quad + \left(C_Z^{2-a} \varep_1^{-(1-a)} W_1^{-(1-a)} W_2^{2-a}  + C_Z^{2-a} W_2\right)u^{\alpha+\lambda(3-2a)-1}.
\end{aligned}
\eeq 

Selecting $\varep_1=c_0/2=2^{a-2}$ and combining \eqref{Ionea}, \eqref{uX1}, \eqref{uX2}, we have 
\begin{align*}
I_1&\ge \frac{c_0}{2}I_0 - C_1\int_U (W_1(x) + W_1^{-(1-a)}(x)W_2^{2-a}(x)+ W_2(x)) u^{\alpha+\lambda(3-2a)-1}\d x\\
&\quad -\frac12\int_U u^{\alpha-\lambda-1}a_N \d x,
\end{align*}
where
\beq\label{Cdef}
C_1=\max\{C_Z^{2-a}, \varep_1^{-(1-a)}C_Z^{2-a}\}=(2^{1-a} C_Z)^{2-a}.
\eeq 
Regarding the second term on the right-hand side of the preceding inequality, we note that 
$$W_1(x) + W_1^{-(1-a)}(x)W_2^{2-a}(x)=W_3(x),$$
and, by Young's inequality with powers $(2-a)$ and $(2-a)/(1-a)$, 
$$W_2=(W_2 W_1^{-(1-a)/(2-a)})\cdot W_1^{(1-a)/(2-a)}\le W_2^{2-a} W_1^{-(1-a)}+W_1=W_3.$$
Therefore, we have
\beq \label{I1est}
I_1\ge \frac{c_0}{2}I_0 - 2 C_1 A_1-\frac12A_2,
\eeq 
where 
$$A_1=\int_U W_3(x) u(x,t)^{\alpha+\lambda(3-2a)-1}\d x,\quad 
A_2=\int_U u(x,t)^{\alpha-\lambda-1}a_N(x) \d x.
$$
Combining \eqref{dJ1} with \eqref{I1est} gives
\beq\label{basicIneq} 
\frac{\lambda}{\alpha}J_0'(t)+ \frac{c_0}{2}(\alpha-\lambda)I_0(t)\le  2C_1(\alpha-\lambda)A_1+\frac12(\alpha-\lambda)A_2 + I_2.
\eeq 
Below, we estimate each term $A_1$, $A_2$ and $I_2$.

\medskip\noindent\textit{Estimate of $A_1$. }
Because $r>\lambda(3-2a)-1$ from \eqref{newrr} and $\alpha>1$ from the first condition in \eqref{newaa}, we have
\beqs
0<r-\lambda(3-2a)+1<r+1<\alpha+r.
\eeqs 
Applying the Young inequality with the powers $\frac{\alpha+r}{r-\lambda(3-2a)+1}$  and $\frac{\alpha+r}{\alpha+\lambda(3-2a)-1}$, we have
\beq\label{A1pre}
2C_1A_1 =\int_U \left(2C_1 W_3(x)\right)\cdot  u^{\alpha+\lambda(3-2a)-1}\d x \le C_2{\mathcal K}_5+ \int_U |u|^{\alpha+r} \d x,   
\eeq
where $C_2=(2C_1)^\frac{\alpha+r}{r-\lambda(3-2a)+1}$.
To estimate the last integral, we apply inequality \eqref{S10} in Lemma \ref{WS1} with the choice \eqref{choice}.
Note, in this case, that  the number $m$ in \eqref{powers_def} becomes 
\beqs 
 m= \frac{\alpha-1+\lambda+a}{2-a},
\eeqs 
and the  numbers $\theta$, $\mu_1$ in \eqref{powers_def} become those in \eqref{mtm}.
Regarding the last condition in \eqref{manyalpha}, we observe that
\beq\label{obs}
\frac{p-s}{p-1}=\frac{1-a-\lambda}{1-a}=1-\frac{\lambda}{1-a}<1.
\eeq  
Thus, condition \eqref{manyalpha} is reduced to $\alpha\ge \lambda+1$ which is met by the first condition in \eqref{newaa}.
Also, condition \eqref{many12} reads as
\beq\label{prepcond}
\alpha>\frac{2(r+a+\lambda-1)}{r_*}.
\eeq
Then, under assumption \eqref{prepcond},  we obtain from \eqref{A1pre} and \eqref{S10}, for any number $\varep_2>0$,  that
\beq\label{A1ineq}
2C_1A_1 \le C_2{\mathcal K}_5 + \varep_2 I_0 + Z_1 (J_0^{1+r/\alpha}+ J_0^{1+\mu_1/\alpha} ),
\eeq
where, referring to the notation in Lemma \ref{WS1} and \eqref{GKs},  
\begin{align*}
Z_1&=\max\left\{
D_{1,m,\theta}\Phi_1,\ 
\varep_2^{-\frac{\theta}{1-\theta}} D_{2,m,\theta} \Phi_2
 \right\},
\end{align*}
with
\beqs
\Phi_1=G_1^\theta G_3^{1-\theta}= {\mathcal K}_2^\frac{\theta(\alpha(1-a)-1+\lambda+a)}{\alpha} {\mathcal K}_1^{(1-\theta)(1+\mu_1/\alpha)} ,\quad 
\Phi_2= G_2^\frac{\theta}{1-\theta} G_3={\mathcal K}_4^{\frac{ \theta(1-r_1)}{r_1(1-\theta)}} {\mathcal K}_1^{1+\mu_1/\alpha}.
\eeqs

\medskip\noindent\textit{Estimate of $A_2$. } 
By H\"older's inequality for the powers $\alpha/(\alpha-\lambda-1)$ and $\alpha/(\lambda+1)$, 
\beq \label{A2ineq}
\begin{aligned}
    A_2&=\int_U \left (u^{\alpha-\lambda-1}\phi^\frac{\alpha-\lambda-1}{\alpha} \right)\cdot \left( \phi^{-\frac{\alpha-\lambda-1}{\alpha}} a_N \right)\d x\\
&\le \left(\int_U u^\alpha\phi \d x\right)^{1-\frac{\lambda +1}{\alpha}} \left (\int_U (a_N\phi^{-\frac{\alpha-\lambda-1}{\alpha}})^{\frac{\alpha}{\lambda+1} } \d x\right )^{\frac{\lambda+1}{\alpha}}= Z_2 J_0^{1-\frac{\lambda +1}{\alpha}},
\end{aligned}
\eeq 
where $Z_2={\mathcal K}_3^{\frac{\lambda+1}{\alpha}}$. 

\medskip\noindent\textit{Estimate of $I_2$. } We apply the Young inequality with the powers $ \frac{\alpha+r}r$ and $ \frac{\alpha+r}\alpha$ to have
\beq\label{I2ineq}
I_2= \int_\Gamma \psi^- u^\alpha \d S \le K+\int_\Gamma u^{\alpha+r}\d S,
\text{ where }
K=\int_\Gamma (\psi^-)^{\frac{\alpha+r}{r}} \d S.
\eeq
To estimate the surface integral of $u^{\alpha+r}$, we apply Lemma \ref{trace} to the same functions $W(x)$, $\varphi(x)$ and numbers $p$, $s$ as in \eqref{choice}.
Again, with observation \eqref{obs}, condition \eqref{manyalpha} is reduced to $\alpha\ge \lambda+1$ which is met by the first condition in \eqref{newaa}.
The number $\widetilde r$ in \eqref{rtilde} becomes that in \eqref{newtr}.
Note from \eqref{newrr} that $r>\lambda(3-2a)-1>\lambda-1$,
hence 
\beqs 
\widetilde r>r+\frac{r+1-\lambda}{1-a}>r>0.
\eeqs 
Thus, we use part \ref{trpos}  of Lemma \ref{trace}.
The numbers  $\widetilde\theta$, $\widetilde \mu_1$ in \eqref{tmtilde} become those in \eqref{trtm2}.
Condition \eqref{only} becomes the second condition in \eqref{newaa}.
If $r+a+\lambda-1\le 0$, then  \eqref{prepcond} holds true thanks to the simple fact $\alpha>0$.  
Otherwise, the second condition in \eqref{newaa} already implies \eqref{prepcond}.
Then we  have from \eqref{trace1}, for any number $\varep_3>0$,  that 
\beq \label{trace30}
\int_\Gamma u^{\alpha+r} \d S 
\le 3\varepsilon_3 I_0+Z_3(J_0^{1+r/\alpha}+J_0^{1+\mu_1/\alpha}+J_0^{1+\widetilde r/\alpha}+J_0^{1+\widetilde \mu_1/\alpha}),
\eeq 
where
\begin{align*}
    Z_3=\max\Big\{& c_5D_{1,m,\theta} \Phi_3,\
    \varep_3^{-\frac{\theta}{1-\theta}}c_5^{\frac 1{1-\theta}}D_{2,m,\theta}\Phi_4,\ 
     \varep_3^{-\frac 1{1-a}}(c_6 (\alpha+r))^\frac {2-a}{1-a} D_{1,m,\widetilde \theta} \Phi_6,\\
     & 
     \varep_3^{-(\frac 1 {1-a}+\frac {2-a}{1-a}\cdot \frac{\widetilde\theta}{1-\widetilde\theta} )}( c_6(\alpha+r))^{\frac {2-a}{1-a}\cdot \frac 1{1-\widetilde\theta}} D_{2,m,\widetilde \theta} \Phi_7
 \Big\},
\end{align*}
with, see \eqref{Ph34}, \eqref{Ph67} and \eqref{GKs}, 
\begin{align*}
\Phi_3&=G_1^\theta G_4^{1-\theta}
    ={\mathcal K}_2^\frac{\theta(\alpha(1-a)-1+\lambda+a)}{\alpha} {\mathcal K}_1^{(1-\theta)(1+\mu_1/\alpha)} ,\ 
\Phi_4= G_2^\frac{\theta}{1-\theta} G_4
    =  {\mathcal K}_4^{\frac{ \theta(1-r_1)}{r_1(1-\theta)}} {\mathcal K}_1^{1+\mu_1/\alpha},\\ 
\Phi_6&=  G_1^{\widetilde\theta} G_5^{1-\widetilde\theta}
= {\mathcal K}_2^\frac{\widetilde\theta(\alpha(1-a)-1+\lambda+a)}{\alpha} {\mathcal K}_6^{(1-\widetilde\theta)(1+\widetilde\mu_1/\alpha)} 
,\ 
\Phi_7=  G_2^{\frac{\widetilde\theta}{1-\widetilde\theta}} G_5= {\mathcal K}_4^{\frac{ \widetilde\theta(1-r_1)}{r_1(1-\widetilde \theta)}} {\mathcal K}_6^{1+\widetilde\mu_1/\alpha}.
\end{align*}

\medskip
Combining inequality \eqref{basicIneq} with estimates \eqref{A1ineq}, \eqref{A2ineq}, \eqref{I2ineq} and \eqref{trace30} gives
\begin{align*}
&\frac{\lambda}{\alpha}J_0'(t)+ \frac{c_0}{2}(\alpha-\lambda)I_0(t)
\le  (\alpha-\lambda)\left [C_2{\mathcal K}_5 + \varep_2 I_0 + Z_1( J_0^{1+r/\alpha} +  J_0^{1+\mu_1/\alpha})\right ]\\
&\quad +\frac12(\alpha-\lambda) Z_2 J_0^{1-\frac{\lambda +1}{\alpha}} + 3\varepsilon_3 I_0
+Z_3(J_0^{1+r/\alpha}+J_0^{1+\mu_1/\alpha}+J_0^{1+\widetilde r/\alpha}+J_0^{1+\widetilde \mu_1/\alpha})+K(t).
\end{align*} 
Let us take $\varep_2 = c_0/8$ and $\varep_3 = c_0(\alpha-\lambda)/24$.  It results in
\begin{align*}
&\frac{\lambda}{\alpha}J_0'(t)+ \frac{c_0}{4}(\alpha-\lambda)I_0(t)
\le  ((\alpha-\lambda)Z_1+Z_3)( J_0^{1+r/\alpha} +  J_0^{1+\mu_1/\alpha})\\
&\quad +\frac12(\alpha-\lambda) Z_2 J_0^{1-\frac{\lambda +1}{\alpha}}
+Z_3(J_0^{1+\widetilde r/\alpha}+J_0^{1+\widetilde \mu_1/\alpha})
+ C_2(\alpha-\lambda){\mathcal K}_5 +K(t).
\end{align*} 

We examine the powers of $J_0$ in realtion with  $\mu_{\min}$ and $\mu_{\max}$ defined in \eqref{mmm2}.
On the one hand, since $\alpha> \lambda+1>0$ and $\mu_1,\widetilde\mu_1>-\alpha$, we have $0<\mu_{\min}<\alpha$. On the other hand, $\mu_{\max}\ge \widetilde r>r>0$.
We then have
\beqs
-\mu_{\min}\le -(\lambda+1),r,\mu_1,\widetilde r,\widetilde \mu_1\le \mu_{\max}.
\eeqs
Then applying inequality \eqref{ee2} yields
$$
J_0^{1-(\lambda+1)/\alpha}, J_0^{1+r/\alpha},J_0^{1+\mu_1/\alpha},
 J_0^{1+\widetilde r/\alpha}, J_0^{1+\widetilde\mu_1/\alpha}
\le J_0^{1-\mu_{\min}/\alpha}+J_0^{1+\mu_{\max}/\alpha} .
$$
Thus,
\beq \label{diff1}
\frac{\lambda}{\alpha}J_0'(t)+ \frac{c_0}{4}(\alpha-\lambda)I_0(t)
\le  Z_4( J_0^{1-\mu_{\min}/\alpha} +  J_0^{1+\mu_{\max}/\alpha})+ C_2(\alpha-\lambda){\mathcal K}_5 +K(t),
\eeq  
where
\beqs
Z_4=2 ((\alpha-\lambda)Z_1+Z_3) + \frac12(\alpha-\lambda)Z_2  +2Z_3.
\eeqs
We deduce from \eqref{diff1} the
desired inequality \eqref{diffIneq} with 
$Z_*=\frac{\alpha}{\lambda} \max\left\{1, Z_4,C_2(\alpha-\lambda){\mathcal K}_5\right\}.$
\end{proof}

\begin{remark}   
The following remarks on Proposition \ref{Diff4u} and its proof are in order.
\begin{enumerate}[label=\rnum]
    \item In the case $\psi\ge 0$, one has $I_2= 0$, and there is no need to deal with the boundary integral. Consequently, the inequality \eqref{diffIneq} will become much simpler.
    \item If $\alpha=\lambda+1$, then
\beqs 
A_2=\int_U a_N(x)\d x={\mathcal K}_3=Z_2 J_0^{1-\frac{\lambda +1}{\alpha}}
\eeqs
with the use of the convention $(J_0)^0=1$.
Thus, we have the same estimate \eqref{A2ineq}.
Hence, in inequality \eqref{diffIneq} above, $\mu_{\min}=\lambda+1=\alpha$, and
$\left( \int_U u^\alpha \phi(x)\d x\right)^{1-{\mu_{\min}}/{\alpha}}$ can be replaced with $1$.
\end{enumerate}
\end{remark}

The main estimates of this section for the solutions of problem \eqref{mainpb} is the following.

\begin{theorem}\label{Labound}
Let the numbers $\alpha$, $\mu_{\max}$
and  $Z_*$ be as in Proposition \ref{Diff4u}, and $M(t)$ be defined by \eqref{Mtdef}. 
Set $V_{0}=1+\int_U \phi(x) u_0(x)^{\alpha}\d x$.
Suppose there is a number $T>0$  such that
 \beq\label{tsmall1}
  \int_0^T M(\tau)\d\tau \le  \frac{\alpha}{3Z_* \mu_{\max} } V_{0}^{-\frac{\mu_{\max}}\alpha} .
\eeq 
\begin{enumerate}[label=\tnum]
    \item Then one has, for all $t\in[0,T)$,
    \beq\label{u_wgt_est}
    \int_U u^\alpha(x,t)\phi(x) \d x \le \left( V_{0}^{-\frac{\mu_{\max}}\alpha}   -\frac{3Z_* \mu_{\max} }\alpha \int_0^t M(\tau) \d\tau  \right)^{-\frac\alpha{\mu_{\max}}}.
    \eeq

\item Consequently, if $\beta$ is a number in the interval $(0,\alpha)$ such that  
$$C_{\alpha,\beta}\eqdef \int_U \phi(x)^{-\frac\beta{\alpha-\beta}} \d x\text{ is finite,}$$
then one has, for all $t\in[0,T)$, 
\beq\label{ube}
    \int_U u^\beta(x,t) \d x \le C_{\alpha,\beta}^{1-\beta/\alpha} \left( V_{0}^{-\frac{\mu_{\max}}\alpha}   -\frac{3Z_* \mu_{\max} }\alpha \int_0^t M(\tau) \d\tau  \right)^{-\frac\beta{\mu_{\max}}}.
    \eeq
\end{enumerate}
\end{theorem}
\begin{proof}
(i) Set $V(t)= 1+\int_U \phi(x) u(x,t)^\alpha \d x$.
Let the number $\mu_{\min}$ be as in Proposition \ref{Diff4u}. Considering inequality \eqref{diffIneq}, we estimate
\begin{align*} 
 \left( \int_U u^\alpha \phi(x)\d x\right)^{1-{\mu_{\min}}/{\alpha}}
 &\le V(t)^{1-{\mu_{\min}}/{\alpha}}\le V(t)^{1+\mu_{\max}/\alpha},\\
 \left( \int_U u^\alpha \phi(x)\d x\right)^{1+\mu_{\max}/\alpha}
 &\le V(t)^{1+\mu_{\max}/\alpha}.
\end{align*} 
Then neglecting the integral $\int_U |u|^{\alpha-\lambda-1}|\nabla u|^{2-a} W_1(x)\d x$ in \eqref{diffIneq}, we obtain
\beqs 
    V'(t) \le Z_*\left (2V(t)^{1+\mu_{\max}/\alpha}+ M(t)\right)
\le Z_*\left (2V(t)^{1+\mu_{\max}/\alpha}M(t) + V(t)^{1+\mu_{\max}/\alpha}M(t)\right),
\eeqs
thus,
\beqs 
    V'(t)  \le 3Z_*M(t)V(t)^{1+\mu_{\max}/\alpha}.
\eeqs 
Solving this differential  inequality gives
\begin{align*}
  V(t) &\le \left (V_{0}^{-\mu_{\max}/\alpha} - \frac{3Z_* \mu_{\max} }\alpha \int_0^t M(\tau)\d \tau\right)^{-\alpha/\mu_{\max}}. 
\end{align*}
Then estimate \eqref{u_wgt_est} follows.

(ii) By H\"older's inequality with powers $\alpha/\beta$ and $\alpha/(\alpha-\beta)$, we have 
\beqs 
\int_U u^\beta  \d x 
= \int_U (u^\beta \phi^{\beta/\alpha} )\cdot \phi^{-\beta/\alpha}   \d x 
\le \left( \int_U u^{\alpha}\phi  \d x \right )^{\beta/\alpha}  \left (\int_U \phi^{-\beta/(\alpha-\beta)} \d x\right)^{1-\beta/\alpha} .
\eeqs
Combining this with the estimate \eqref{u_wgt_est}, we obtain \eqref{ube}.
\end{proof}

\begin{example}\label{ex1}
For the two-term Forchheimer law \eqref{2term}, we have $N=1$, $\bar\alpha_0=0$ and $\bar\alpha_N=1$ in equation \eqref{Fhetero}. Then the number $a$ in \eqref{eq9} is $a=1/2$.
In the two and three dimensional cases, the number $r_1$ in \eqref{newro} and $r_*$ in \eqref{newrs} satisfy
\beq\label{rr23}
\begin{aligned}
    &\frac23\le r_1<1\text{ and } r_*=\frac74-\frac1{r_1}\text{ for }n=2,\\
    &\frac23< r_1<1 \text{ and } r_*=\frac32-\frac1{r_1} \text{ for }n=3.
\end{aligned}
\eeq

Now, consider the ideal gas. Then $\gamma=1$ and, by  \eqref{ulam}, $\lambda=1/2$. 
Note in this case that $1-\lambda-a=0$ and $\lambda(3-2a)-1=0$. 
From \eqref{newrr} and \eqref{newtr},
\beqs
r>0\text{ and } \widetilde r=3r. 
\eeqs
Condition \eqref{newaa} becomes
\beqs
\alpha>\max\left\{\frac32,\frac{6r}{r_*}\right\}.
\eeqs
From \eqref{mtm} and \eqref{trtm2},
\beq\label{TM}
\theta=\Theta(r)=\frac{\alpha+2r}{\alpha(1+r_*)},\quad \mu_1=\Lambda(r,\theta)=\frac{r}{1-\theta},
\eeq
\beq \label{TM1}
\widetilde\theta=\Theta(\widetilde r)=\frac{\alpha+6r}{\alpha(1+r_*)},\quad 
\widetilde\mu_1=\Lambda(\widetilde r,\widetilde \theta)=\frac{3r}{1-\widetilde\theta} =\frac{3r}{1-\frac{\alpha+6r}{\alpha(1+r_*)}}
=\frac{3r(1+r_*)\alpha}{r_*\alpha-6r}.
\eeq 

Then $\mathcal K_1,\ldots,\mathcal K_6$ in \eqref{KKs} become
\begin{align*}
{\mathcal K}_1&= \int_U  \phi(x)^{-1}\d x,\quad 
{\mathcal K}_2= \int_U \phi(x)^{-2} \d x,\quad 
{\mathcal K}_3= \int_U a_N(x)^{\frac{2\alpha}{3} }\phi(x)^{1-\frac{2\alpha}{3} } \d x\\
{\mathcal K}_4&= \int_U W_1(x)^{-\frac{r_1}{1-r_1}}\d x,\quad 
{\mathcal K}_5= \int_U W_3(x)^{\frac{\alpha+r}{r}} \d x,\\
 {\mathcal K}_6&= \int_U  \phi(x)^{-1}(x)W_1(x)^{-2\alpha/[(1-\widetilde\theta)\alpha+\widetilde r]}\d x
 =\int_U  \phi(x)^{-1}(x)W_1(x)^{-\frac{2(1+r_*)\alpha}{r_*\alpha-3r(1-r_*)}}\d x.
\end{align*}
Clearly, $\widetilde r>r$, and by the particular formula of $\Theta(\cdot)$ in \eqref{TM},
$\widetilde\theta=\Theta(\widetilde r)>\Theta(r)=\theta$.
Together the fact that $\Lambda(r,\theta)$ in \eqref{TM} is strictly increasing in $r>0$ and $\theta\in(0,1)$, these imply
\beqs
 \widetilde\mu_1=\Lambda(\widetilde r,\widetilde \theta)>\Lambda(r,\theta)=\mu_1.
\eeqs
Note also that 
\beqs
\widetilde\mu_1=\Lambda(\widetilde r,\widetilde \theta)=\frac{\widetilde r}{1-\widetilde\theta}>\widetilde r.
\eeqs
Therefore, the number $\mu_{\max}$ in \eqref{mmm2} simply is $\mu_{\max}=\widetilde \mu_1$. Using this relation and \eqref{TM1} gives
\beq\label{mumax}
\mu_{\max}=\frac{3r(1+r_*)\alpha}{r_*\alpha-6r}.
\eeq
Thus, the last power in estimate \eqref{u_wgt_est}  is
\beqs
-\frac {\alpha}{\mu_{\max}} = -\frac {r_*\alpha-6r}{3r(1+r_*)}.
\eeqs 
\end{example}

\section{Estimates for the space-time essential supremum}\label{maxsec}

In this section, we use Moser's iteration to estimate the $L_{x,t}^\infty$-norm of the solution $u(x,t)$ of the problem \eqref{mainpb}. Throughout this section, $u(x,t)$ is a non-negative solution of the problem \eqref{mainpb} as in Section \ref{lasec}.

Below, we denote the positive constants  $\widetilde c_0=2^{a-2}$ and $\widetilde C_1=2C_1$, see \eqref{Cdef}. 
Recall that $c_5$, $c_6$ are the positive constants from the trace theorem \eqref{firstrace} and  the weight $W_3(x)$ is defined by \eqref{W3}. 
For $T>0$, denote $Q_T =U\times(0,T)$.

\begin{lemma} \label{caccio} 
Let numbers $\widetilde \kappa>1$ and $p_i>1$ for $i=1,2,\ldots,5$  satisfy
\beq\label{ktil}
p_1,p_2,p_3p_4<\widetilde \kappa\text{ and } p_6\eqdef \frac{p_5[p_3(2-a)-1]}{1-a}<\widetilde \kappa.
\eeq 
For $i=1,2\ldots, 5$, let $q_i$  be the H\"older conjugate exponent of $p_i$. 
Assume the integrals 
$\int_U\phi(x)\d x$ and 
\beq\label{wcond1}
\int_U\left\{\phi(x)^{-\frac{q_4}{p_4}}+ W_3(x)^{q_1}\phi(x)^{-\frac{q_1}{p_1}}+a_N(x)^{q_2}\phi(x)^{-\frac{q_2}{p_2}}+ W_1(x)^{-\frac{q_5}{1-a}}\phi(x)^{-\frac{q_5}{p_5}}\right\} \d x
\eeq 
are finite.
Let $T>T_2>T_1\ge 0$.
If $\alpha$ is a number satisfying 
\beq\label{alp-Large}
\alpha >  \max\left\{\lambda+1, \frac{p_1(\lambda(3-2a)-1)}{\widetilde \kappa-p_1}, \frac{p_5(a+\lambda-1)}{(1-a)(\widetilde\kappa -p_6)}\right\},
\eeq
then one has  
\beq\label{newS1}
\begin{aligned}
&\essup_{t\in(T_2,T)} \int_U u^\alpha(x,t)\phi(x)\d x \\
&\le  c_8\alpha^{\frac{3-2a}{1-a}}(1+T)\left(1+\frac1{T_2-T_1}\right) \mathcal N_1 \Psi_T 
\Big(\norm {u}_{L_\phi^{\widetilde \kappa\alpha}(U\times(T_1,T))}^{\alpha-h_1}+\norm {u}_{L_\phi^{\widetilde \kappa\alpha}(U\times(T_1,T))}^{\alpha+h_2}\Big),
\end{aligned}
\eeq
and
\beq\label{newS2}
\begin{aligned}
&\int_{T_2}^T\int_U u^{\alpha-\lambda-1}(x,t)|\nabla u(x,t)|^{2-a}W_1(x)\d x \d t\\
&\le c_9\alpha^{\frac{2-a}{1-a}}(1+T)\left(1+\frac1{T_2-T_1}\right) \mathcal N_1 \Psi_T  
\Big(\norm {u}_{L_\phi^{\widetilde \kappa\alpha}(U\times(T_1,T))}^{\alpha-h_1}+\norm {u}_{L_\phi^{\widetilde \kappa\alpha}(U\times(T_1,T))}^{\alpha+h_2}\Big),
\end{aligned}
\eeq
where $\mathcal N_1$ and $\Psi_T$ are defined by
\begin{align}
\mathcal N_1
&=\left(1+\int_U\phi(x) \d x\right)\Bigg\{ 1+ \left(\int_U W_3(x)^{q_1}\phi(x)^{-\frac{q_1}{p_1}} \d x\right)^\frac{1}{q_1}+\left(\int_U a_N(x)^{q_2}\phi(x)^{-\frac{q_2}{p_2}} \d x\right)^\frac{1}{q_2} \notag \\
&\quad  +\left(\int_U \phi(x)^{-\frac{q_4}{p_4}} \d x\right)^\frac{1}{p_3q_4}
+\left(\int_U W_1(x)^{-\frac{q_5}{1-a}}\phi(x)^{-\frac{q_5}{p_5}} \d x\right)^\frac{1-a}{q_5(p_3(2-a)-1)}
\Bigg\}, \label{N1} \\
\label{bcquant}
\Psi_T&= 1+\left(\int_0^T \int_\Gamma (\psi^-(x,t))^{q_3} \d S \d t\right)^{\frac{p_3(2-a)}{q_3[p_3(2-a)-1]}},
\end{align} 
the numbers $h_1>1$ and $h_2\ge 0$ are
\beq \label{h12}
h_1=\lambda+1, \quad 
h_2=\max\left\{0,3(\lambda-2a)-1,\frac{a+\lambda-1}{p_3(2-a)-1}\right\},
\eeq 
and the constants $c_8>1$ and $c_9>0$ are 
\beq\label{c67}
c_8=\lambda^{-1}c_{10},\ 
c_9=8c_{10} \text{ with }
c_{10} = 20\max\left\{ \widetilde C_1, 4\lambda,c_5^{1/p_3}, \left[4(2c_6p_3)^{2-a}\right]^{\frac 1{p_3(2-a)-1}}\right\}.
\eeq 
\end{lemma}
\begin{proof}
It is obvious from \eqref{ktil} that $p_3,p_4<\widetilde \kappa$ and $1<p_5<p_6<\widetilde \kappa$.
Let $\xi=\xi(t)$ be a $C^1$-function  on $[0,T]$ with $0\le \xi(t)\le 1$ on $[0,T]$, and  
\beq\label{xiprop} 
\xi(t)=0 \text{ on } [0,T_1], \ 
\xi(t)=1 \text{ on } [T_2,T], \text{ and }
 0\le \xi'(t)\le \frac{2}{T_2-T_1} \text{ on } [0,T].
\eeq 

Multiply the PDE in \eqref{mainpb} by test function $u^{\alpha-\lambda}\xi^2$ and  integrating the resulting equation over $U$ give
 \begin{align*}
\frac{\lambda}{\alpha}\ddt \int_U \phi(x) u^{\alpha} \xi^2 \d x
-\frac{\lambda}{\alpha}\int_U 2\phi(x)u^{\alpha} \xi \xi' \d x
= \int_U \nabla\cdot \mathbf X(x,\nabla u+\mathbf Z(u)) u^{\alpha-\lambda}\xi^2 \d x.
\end{align*}
Performing integration by parts for the last integral and using the boundary condition give
 \begin{align*}
\frac{\lambda}{\alpha}\ddt \int_U \phi u^{\alpha} \xi^2 \d x
&=-(\alpha-\lambda) \int_U \mathbf X(x,\nabla u+\mathbf Z(u)) \cdot (\nabla u)u^{\alpha-\lambda-1}\xi^2 \d x\\
&\quad - \int_\Gamma \psi u^\alpha \xi^2 \d x
+\frac{2\lambda}{\alpha}\int_U \phi(x)u^{\alpha} \xi \xi' \d x.
\end{align*}

For the first integral on the right-hand side using the similar estimate to \eqref{I1est}, we obtain 
\beq \label{bint}
\begin{aligned}
&\frac{\lambda}{\alpha}\ddt\int_U \phi u^\alpha \xi^2 \d x
\le  - \widetilde c_0(\alpha-\lambda)\int_U W_1(x)  u^{\alpha-\lambda-1} |\nabla u|^{2-a} \xi^2 \d x \\
&\quad + \widetilde C_1(\alpha-\lambda)\int_U W_3(x) u^{\alpha+\lambda(3-2a)-1}\xi^2 \d x
 +\frac12(\alpha-\lambda)\int_U u^{\alpha-\lambda-1}a_N \xi^2 \d x \\
&\quad +\frac{2\lambda}{\alpha}\int_U \phi(x)u^{\alpha} \xi \xi' \d x
+ \int_\Gamma \psi^- u^\alpha \xi^2 \d S.
\end{aligned}
\eeq 
Denote
\beqs 
B_0=\essup_{t\in(0,T)}\int_U \phi u^\alpha \xi^2 \d x,\quad 
A_0=\iint_{Q_T} W_1(x)  u^{\alpha-\lambda-1} |\nabla u|^{2-a} \xi^2 \d x\d t.
\eeqs 
Note that
\beqs
\int_U \phi(x) u^\alpha(x,0) \xi^2(0) \d x=0.
\eeqs
Then integrating inequality \eqref{bint} in time from $0$ to $t$, for any $t\in(0,T)$, and  taking  the essential supremum with respect to $t$ yield
\beq\label{Bzest}
\frac{\lambda}{\alpha}B_0
\le \mathcal I_0\eqdef   \widetilde C_1(\alpha-\lambda)A_1 + \frac12(\alpha-\lambda)A_2+\frac{2\lambda}{\alpha}A_3 + A_4,
\eeq
where 
\begin{align*}
    A_1&=\iint_{Q_T} W_3(x) u^{\alpha+\lambda(3-2a)-1}\xi^2\d x\d t,&& 
    A_2=\iint_{Q_T} u^{\alpha-\lambda-1}a_N\xi^2 \d x \d t,\\
A_3&=\iint_{Q_T} \phi(x) u^{\alpha} \xi \xi' \d x \d t,&&
A_4=\int_0^T \int_{\Gamma} \psi^- u^\alpha \xi^2 \d S \d t.
\end{align*}
Now, just integrating inequality \eqref{bint} in time from $0$ to $T$ gives  
\beq \label{Azest}
\widetilde c_0(\alpha-\lambda) A_0\le \mathcal I_0.
\eeq
Adding \eqref{Bzest} to \eqref{Azest} yields
\beq\label{MainD}
\frac{\lambda}{\alpha}B_0
+\widetilde c_0(\alpha-\lambda) A_0
\le 2\mathcal I_0=   2\widetilde C_1(\alpha-\lambda)A_1 + (\alpha-\lambda)A_2+\frac{4\lambda}{\alpha}A_3 + 2A_4.
\eeq

Denote $E(\alpha)=\iint_{Q_T} \phi(x)u^{\alpha}(x,t)\xi(t)  \d x \d t$.  We will use the fact $0\le \xi^2\le \xi\le 1$.

\medskip
\noindent\textit{Estimation of $A_1$ and $A_2$.}
Let 
\beqs 
\alpha_1=p_1(\alpha+\lambda(3-2a)-1)\text{ and } \alpha_2=p_2(\alpha-\lambda-1).
\eeqs 
By applying H\"older's inequality with the powers $p_1$ and $q_1$, we have
\beq\label{A1est}
\begin{aligned}
A_1&=\iint_{Q_T}\left( u^{\alpha+\lambda(3-2a)-1}\phi^{1/p_1}\right)\cdot \left(\phi^{-1/p_1}  W_3(x)  \right) \xi^2\d x\d t\\
&\le \left(\iint_{Q_T} \phi(x) u^{\alpha_1}\xi^2 \d x \d t\right)^{1/p_1} F_1 \le E(\alpha_1)^{1/p_1} F_1,    
\end{aligned}
\eeq
where 
\beqs
	F_1=\left(\iint_{Q_T} W_3^{q_1} \phi^{-q_1/p_1}\xi^2 \d x\d t\right)^{1/q_1}.
\eeqs 
Similarly, by applying H\"older's inequality with the powers $p_2$ and $q_2$, we have
\beq\label{A2est}
\begin{aligned}
A_2& = \iint_{Q_T} \left(\phi^{1/p_2}u^{\alpha-\lambda-1}\right)\cdot \left(\phi^{-1/p_2} a_N\right) \xi^2 \d x \d t  \\
&\le \left(\iint_{Q_T} \phi(x) u^{\alpha_2}\xi^2 \d x \d t\right)^{1/p_2} F_2 \le E(\alpha_2)^{1/p_2}F_2,    
\end{aligned}
\eeq
where 
\beqs
	F_2=\left(\iint_{Q_T} a_N^{q_2}\phi^{-q_2/p_2} \xi^2\d x\d t\right)^{1/q_2}.
\eeqs 

\medskip
\noindent\textit{Estimation of $A_3$.}
By \eqref{xiprop},  we have 
\beq\label{IIx}
A_3=\iint_{Q_T} \phi(x) u^{\alpha} \xi \xi' \d x \d t
\le \frac{2}{T_2-T_1}\iint_{Q_T} \phi(x) u^{\alpha} \xi  \d x \d t =\frac{2}{T_2-T_1} E(\alpha).
\eeq

\medskip
\noindent\textit{Estimation of $A_4$.}
Let $\alpha_3=p_3 \alpha$. We first apply H\"older's inequality with the powers $p_3$ and $q_3$ and obtain
\beq\label{Afour}
    A_4
        \le \left(\int_0^T \int_{\Gamma} u^{p_3\alpha} \xi^2\d S\d t\right)^{1/p_3} \left(\int_0^T \int_{\Gamma}  (\psi^-)^{q_3} \xi^2\d S\d t\right)^{1/q_3} 
    =J^{1/p_3} F_3 ,
\eeq 
 where 
 \beqs 
 J= \int_0^T \int_{\Gamma} u^{\alpha_3} \xi^2\d S\d t, \quad F_3= \left(\int_0^T \int_{\Gamma} (\psi^-)^{q_3} \xi^2\d S\d t\right)^{1/q_3}.
 \eeqs 
Now applying the trace theorem \eqref{firstrace} to the power $p:=1$ and function $f:=u^{\alpha_3}$, we have
\begin{equation}\label{i4est}
    J \le c_5 \iint_{Q_T} u^{\alpha_3}\xi^2 \d x \d t +c_6\alpha_3 \iint_{Q_T} u^{\alpha_3-1}|\nabla u|\xi^2 \d x \d t.
\end{equation}
For the first integral on the right-hand side of \eqref{i4est}, letting 
\beqs 
\alpha_4=p_4\alpha_3=p_4p_3\alpha
\eeqs 
and applying H\"older's inequality with the powers $p_4$ and $q_4$,
we have 
\beq\label{I4part1}
 \iint_{Q_T} u^{\alpha_3}\xi^2 \d x \d t
= \iint_{Q_T} (u^{\alpha_3}\phi^{1/p_4})\cdot \phi^{-1/p_4}\xi^2 \d x \d t 
\le E(\alpha_4)^{1/p_4} F_4,
\eeq
where
 \beqs 
 F_4=\left(\iint_{Q_T} \phi^{-q_4/p_4}\xi^2 \d x \d t\right)^{1/q_4}.
 \eeqs 
To  estimate the second integral on the right-hand side of \eqref{i4est}, we set the numbers
\begin{align}
m_4&=\frac{(\alpha_3-1)(2-a)-(\alpha-\lambda-1)}{1-a}
=\alpha \frac{p_3(2-a)-1}{1-a} +\frac{a+\lambda-1}{1-a},\notag \\
\label{al5}
\alpha_5&=p_5 m_4=\alpha p_6 + \frac{p_5(a+\lambda-1)}{1-a}.    
\end{align}
Applying H\"older's inequality with the powers $(2-a)$ and $(2-a)/(1-a)$, we have 
\begin{align*}
&\iint_{Q_T} u^{\alpha_3-1}|\nabla u|\xi^2 \d x \d t
=\iint_{Q_T}\left( u^\frac{\alpha-\lambda-1}{2-a} |\nabla u| W_1^\frac1{2-a}\right)
\cdot \left(  u^{\alpha_3-1-\frac{\alpha-\lambda-1}{2-a}} W_1^{-\frac1{2-a}}\right)   \xi^2 \d x \d t\\
&
 \le  \left(\iint_{Q_T} W_1(x)  u^{\alpha-\lambda-1} |\nabla u|^{2-a} \xi^2 \d x \d t\right)^{1/(2-a)}
    \left(\iint_{Q_T}W_1^{-\frac1{1-a}} u^{m_4}\xi^2 \d x \d t\right)^{(1-a)/(2-a)}\\
    &= A_0^{1/(2-a)} \left(\iint_{Q_T} 
    (\phi^\frac{1}{p_5}u^{m_4})\cdot (W_1^{-\frac1{1-a}}\phi^{-\frac{1}{p_5}}) \xi^2 \d x \d t\right)^{(1-a)/(2-a)}. 
\end{align*}
Applying H\"older's inequality again with the powers $p_5$ and $q_5$ to estimate the last integral, we obtain 
\beq\label{I4part2}
\begin{aligned}
&\iint_{Q_T} u^{\alpha_3-1}|\nabla u|\xi^2 \d x \d t\\
    &\le A_0^{1/(2-a)}
    \Big(\iint_{Q_T}\phi u^{\alpha_5}\xi^2 \d x \d t\Big)^{\frac{1-a}{p_5(2-a)}}
    \Big(\iint_{Q_T}W_1^{-\frac{q_5}{1-a}}\phi^{-q_5/p_5}\xi^2 \d x \d t\Big)^{\frac{1-a}{q_5(2-a)}}\\
    &    \le A_0^{1/(2-a)}E(\alpha_5)^{\frac{1-a}{p_5(2-a)}} F_5^{\frac{1-a}{2-a}}, 
\end{aligned}
\eeq
where
\beqs 
F_5=\left(\iint_{Q_T}W_1^{-\frac{q_5}{1-a}}\phi^{-q_5/p_5}\xi^2 \d x \d t\right)^\frac{1}{q_5}.
\eeqs 
Combining \eqref{i4est} with \eqref{I4part1}, \eqref{I4part2}  and then raising to the power $1/p_3$, we obtain, with the aid of inequality \eqref{ee2} for $p:=1/p_3$,
\begin{equation}\label{i4estf}
    J^{1/p_3} \le c_5^{1/p_3} E(\alpha_4)^{1/p_3p_4} F_4^{1/p_3}+(c_6\alpha_3 )^{1/p_3}A_0^{1/p_3(2-a)}E(\alpha_5)^{\frac{1-a}{p_3p_5(2-a)}} F_5^{\frac{1-a}{p_3(2-a)}}.
\end{equation}
Now using estimate \eqref{i4estf} in \eqref{Afour}, and applying Young's inequality to the last resulting term with powers $p_3(2-a)$ and $\displaystyle\frac{p_3(2-a)}{p_3(2-a)-1}$, we have 
\begin{align*}
    &A_4 \le c_5^{1/p_3} E(\alpha_4)^\frac{1}{p_3p_4}
    F_4^{1/p_3}F_3
    + (c_6p_3\alpha)^{1/p_3} A_0^\frac{1}{p_3(2-a)}
    E(\alpha_5)^{\frac{1-a}{p_3p_5(2-a)}}F_5^{\frac{1-a}{p_3(2-a)}} F_3\\
    & = c_5^{1/p_3} E(\alpha_4)^\frac{1}{p_3p_4}F_6
    +  \left(\varep^\frac{1}{p_3(2-a)} A_0^\frac{1}{p_3(2-a)}\right)\cdot 
    \left( \varep^{-\frac{1}{p_3(2-a)}} (c_6p_3\alpha)^{1/p_3} 
    E(\alpha_5)^{\frac{1-a}{p_3p_5(2-a)}}F_5^{\frac{1-a}{p_3(2-a)}} F_3\right)\\
    &\le c_5^{1/p_3} E(\alpha_4)^\frac{1}{p_3p_4}F_6 
    + \varep A_0
    +\varep^{-\frac{1}{p_3(2-a)-1}}
\left\{  (c_6p_3\alpha)^{1/p_3}  E(\alpha_5)^{\frac{1-a}{p_3p_5(2-a)}} F_5^{\frac{1-a}{p_3(2-a)}} F_3 \right\}^\frac{p_3(2-a)}{p_3(2-a)-1},
\end{align*}
thus,
\beq\label{A4est} 
A_4\le c_5^{1/p_3} E(\alpha_4)^\frac{1}{p_3p_4}F_6 
    + \varep A_0 + \varep^{-\frac{1}{p_3(2-a)-1}}  (c_6 p_3\alpha)^{\frac{2-a}{p_3(2-a)-1}}E(\alpha_5)^{1/p_6}F_7,
\eeq 
where 
\beq\label{F7andF8}
F_6= F_4^{1/p_3}F_3 \text { and } F_7=F_5^\frac{p_5}{p_6} F_3^\frac{p_3(2-a)}{p_3(2-a)-1}.
\eeq

Using \eqref{A1est}, \eqref{A2est} \eqref{IIx} and \eqref{A4est} in \eqref{MainD}, and selecting $\varep=\widetilde c_0(\alpha-\lambda)/4 $, we obtain
\begin{multline}\label{Imain5}
\frac{\lambda}{\alpha}B_0
+ \widetilde c_0\frac{\alpha-\lambda}2A_0
\le 2\widetilde C_1(\alpha-\lambda)E(\alpha_1)^{1/p_1}F_1+\widetilde C_1(\alpha-\lambda)E(\alpha_2)^{1/p_2}F_2+ \frac{8\lambda}{\alpha(T_2-T_1)}E(\alpha)\\
+2c_5^{1/p_3} E(\alpha_4)^{1/(p_3p_4)}F_6+2 \Big(\frac4{\widetilde c_0(\alpha-\lambda)}\Big)^{\frac 1{p_3(2-a)-1}} (c_6  p_3\alpha)^{\frac {2-a}{p_3(2-a)-1}}E(\alpha_5)^{1/p_6}F_7.
\end{multline}

We bound the left-hand side of \eqref{Imain5} from below by 
\begin{align*}
 B_0&\ge \essup_{t\in(T_2,T)}\int_U \phi(x) u^{\alpha}(x,t)\xi^2(t) \d x
=B_*\eqdef \essup_{t\in(T_2,T)}\int_U \phi(x) u^{\alpha}(x,t)\d x,\\
A_0&\ge \int_{T_2}^T\int_U W_1(x)|\nabla u(x,t)|^{2-a}u(x,t)^{\alpha-\lambda-1} \xi^2(t) \d x \d t\\
&=A_*\eqdef \int_{T_2}^T\int_U W_1(x)|\nabla u(x,t)|^{2-a}u(x,t)^{\alpha-\lambda-1} \d x \d t.
 \end{align*}
Thus,
\beq\label{leftest}
\frac{\lambda}{\alpha}B_0+\frac{\widetilde c_0(\alpha-\lambda)}2 A_0
\ge \mathcal I_*\eqdef \frac{\lambda}{\alpha}B_*+\frac{\widetilde c_0(\alpha-\lambda)}2 A_*
 \eeq

We now bound the right-hand side of \eqref{Imain5} from above. We will estimate many $E(\cdot)$-terms by using
\beqs
J_0\eqdef \norm {u}_{L_\phi^{\widetilde \kappa\alpha}(U\times(T_1,T))}=\left(\int_{T_1}^T\int_U \phi u^{\widetilde \kappa \alpha}\d x\d t\right)^{\frac{1}{\widetilde\kappa \alpha}}.
\eeqs
If $0<\beta<\widetilde \kappa \alpha$, then, by H\"older's inequality with powers $\widetilde\kappa \alpha/\beta$ and its H\"older conjugate, we have
\beq\label{X}
E(\beta)\le \int_{T_1}^T \int_U \phi u^\beta  \d x \d t \le J_0^\beta \left(\int_{T_1}^T \int_U \phi \d x\d t\right)^{1-\frac{\beta}{\widetilde \kappa \alpha}} \le  J_0^\beta (T+1) \Phi_*,
\eeq
where
\beq\label{Phistar}
\Phi_*=1+\int_U\phi \d x.
\eeq
From  \eqref{alp-Large}, we have $\alpha>\frac{p_1[\lambda(3-2a)-1]}{\widetilde\kappa-p_1}$, hence
$\alpha_1=p_1(\alpha+\lambda(3-2a)-1)<\widetilde\kappa\alpha$.
From \eqref{ktil}, we have $p_2<\widetilde\kappa$, hence $\alpha_2=p_2(\alpha-\lambda-1)<p_2\alpha<\widetilde\kappa\alpha$.
From \eqref{ktil}, we have $p_3p_4<\widetilde\kappa$, hence $\alpha_4=p_3p_4\alpha<\widetilde\kappa\alpha$.
From  \eqref{alp-Large}, we have $\alpha>\frac{p_5(a+\lambda-1)}{(1-a)(\widetilde\kappa-p_6)}$, hence
$\alpha_5=\alpha p_6 + \frac{p_5(a+\lambda-1)}{1-a}<\widetilde\kappa\alpha$.
Thus, $0<\alpha_1,\alpha_2, \alpha_4,\alpha_5 < \widetilde \kappa\alpha$.
Applying inequality \eqref{X} to $\beta=\alpha$ and $\beta=\alpha_i$, for $i=1,2,4,5$, gives
\beq\label{Eaa}
E(\alpha)\le J_0^\alpha(T+1)\Phi_*\text{ and } E(\alpha_i)\le J_0^{\alpha_i}(T+1)\Phi_*\text{ for }i=1,2,4,5.
\eeq

Next, we estimate
\beq\label{FFT}
F_i\le T^{1/q_i}\widetilde F_i \text{ for }i=1,2,4,5, \text{ and }
F_3\le \widetilde F_3,
\eeq
where 
\begin{align*}
\widetilde   F_1&=\left(\int_U W_3^{q_1} \phi^{-q_1/p_1}\d x\right)^{1/q_1},\ 
\widetilde F_2=\left(\int_U a_N^{q_2}\phi^{-q_2/p_2} \d x\right)^{1/q_2},\\
\widetilde  F_3= \left(\int_0^T \int_{\Gamma} (\psi^-)^{q_3} \d S\d t\right)^{1/q_3},\
\widetilde  F_4&=\left(\int_U \phi^{-q_4/p_4}\d x\right)^{1/q_4},\ 
\widetilde F_5=\left(\int_UW_1^{-\frac{q_5}{1-a}}\phi^{-q_5/p_5}\d x\right)^{\frac{1}{q_5}}.
\end{align*}
Consequently, by combining \eqref{FFT} with the definition \eqref{F7andF8} of $F_6$ and $F_7$, we have 
\beq\label{F78ineq}
F_6\le T^{1/(p_3q_4)}\widetilde  F_4^{1/p_3}\widetilde  F_3,\quad 
F_7\le T^\frac{p_5}{q_5 p_6} \widetilde F_5^\frac{p_5}{p_6} \widetilde F_3^\frac{p_3(2-a)}{p_3(2-a)-1}
=T^\frac{p_5-1}{p_6} \widetilde F_5^\frac{p_5}{p_6} \widetilde F_3^\frac{p_3(2-a)}{p_3(2-a)-1}.
\eeq

For the constants on the right-hand side of \eqref{Imain5}, we denote
\beq\label{Cbar}
C_3 = \max\left\{ 2\widetilde C_1, 8\lambda,2c_5^{1/p_3}, 2\left(\frac {4(c_6 p_3)^{2-a}}{\widetilde c_0}\right)^{\frac 1{p_3(2-a)-1}}\right\}.
\eeq 

Then we have from \eqref{Imain5}, \eqref{leftest}, \eqref{Eaa},  \eqref{FFT}, \eqref{F78ineq} and \eqref{Cbar} that
\begin{multline}\label{Is1}
\mathcal I_*
 \le  C_3\Bigg\{(\alpha-\lambda) J_0^{\alpha_1/p_1}(T+1)^{1/p_1}\Phi_*^{1/p_1} T^{1/q_1}\widetilde F_1
 + (\alpha-\lambda) J_0^{\alpha_2/p_2}(T+1)^{1/p_2} \Phi_*^{1/p_2}T^{1/q_2} \widetilde F_2\\
 + \frac{1}{\alpha(T_2-T_1)}J_0^{\alpha}(T+1)\Phi_*
+J_0^{\alpha_4/(p_3p_4)}(T+1)^{1/(p_3p_4)}\Phi_*^{1/(p_3 p_4)}T^{1/(p_3q_4)}\widetilde  F_4^{1/p_3}\widetilde  F_3\\
 +\Big(\frac{\alpha^{2-a}}{\alpha-\lambda}\Big)^{\frac 1{p_3(2-a)-1}}J_0^{\frac{\alpha_5}{p_6}}(T+1)^{1/p_6}\Phi_*^{1/p_6} T^\frac{p_5-1}{p_6} \widetilde F_5^\frac{p_5}{p_6} \widetilde F_3^\frac{p_3(2-a)}{p_3(2-a)-1}\Bigg\}.
 \end{multline}
Observe, for $i=1,2$, that
\beqs
(T+1)^{1/p_i}T^{1/q_i}\Phi_*^{1/p_i} \le (T+1)^{1/p_i+1/p_i}\Phi_*=(T+1)\Phi_*.
\eeqs
Similarly, 
\beqs  
(T+1)^{1/(p_3p_4)}T^{1/(p_3q_4)}\Phi_*^{1/(p_3 p_4)}\le (T+1)^{1/p_3} \Phi_*^{1/p_3}\le (T+1)\Phi_*.
\eeqs 
Also, thanks to the fact $p_5<p_6$,
\beqs 
(T+1)^{1/p_6}T^\frac{p_5-1}{p_6}\Phi_*^{1/p_6}  \le (T+1)^{p_5/p_6}\Phi_*\le (T+1)\Phi_*.
\eeqs 
Moreover, we estimate, for the first two terms, respectively, the third term, on the right-hand side of \eqref{Is1}
\beqs
\alpha-\lambda<\alpha<\alpha^\frac{2-a}{1-a},\text{ respectively, }
\frac1{\alpha(T_2-T_1)}\le 1+\frac1{T_2-T_1}.
\eeqs
For the last term on the right-hand side of \eqref{Is1}, we estimate, thanks to the facts $\alpha-\lambda\ge 1 $ and $p_3>1$,
\beqs
\Big(\frac{\alpha^{2-a}}{\alpha-\lambda}\Big)^{\frac 1{p_3(2-a)-1}}
\le \alpha^\frac{2-a}{(2-a)-1}=\alpha^\frac{2-a}{1-a}.
\eeqs
Therefore, 
\begin{align*} \mathcal I_*
&\le C_3\alpha^{\frac{2-a}{1-a}}(T+1)\Big(1+\frac {1}{T_2-T_1}\Big)\Phi_* \\
&\quad\times 
\Big (J_0^\frac{\alpha_1}{p_1}\widetilde F_1
+J_0^\frac{\alpha_2}{p_2}\widetilde F_2
+J_0^{\alpha}  + J_0^{\alpha_4/(p_3p_4)}  \widetilde  F_4^{1/p_3}\widetilde  F_3+J_0^{\frac{\alpha_5}{p_6}}\widetilde F_5^\frac{p_5}{p_6} \widetilde F_3^\frac{p_3(2-a)}{p_3(2-a)-1}
\Big)
\end{align*}
which implies
\beq\label{I0fm}
\mathcal I_* \le C_3 \alpha^{\frac{2-a}{1-a}}(T+1)\Big(1+\frac {1}{T_2-T_1}\Big)\Phi_*\widetilde F\cdot\widetilde J,
 \eeq
 where
\begin{align*}
\widetilde F &= \widetilde F_1 +\widetilde F_2+1+\widetilde  F_4^{1/p_3}\widetilde  F_3+\widetilde F_5^\frac{p_5}{p_6} \widetilde F_3^\frac{p_3(2-a)}{p_3(2-a)-1}, \\
\widetilde J &=J_0^{\alpha_1/p_1}+J_0^{\alpha_2/p_2} + J_0^{\alpha} 
+J_0^{\alpha_4/(p_3p_4)}+J_0^{\frac{\alpha_5}{p_6}}.
\end{align*}

Separating the the weights and the boundary data, we estimate $\widetilde F$ by 
\beqs
\widetilde F 
\le \left (1+\widetilde F_1 +\widetilde F_2+\widetilde  F_4^{1/p_3}+\widetilde F_5^\frac{p_5}{p_6}\right)
\left(1+\widetilde  F_3+\widetilde F_3^\frac{p_3(2-a)}{p_3(2-a)-1}\right).
\eeqs
Simply estimating, by \eqref{ee5}, $\widetilde  F_3\le 1+\widetilde F_3^\frac{p_3(2-a)}{p_3(2-a)-1}=\Psi_T$, we obtain 
\beq\label{tildeFfm}
\widetilde F 
\le \left (1+\widetilde F_1 +\widetilde F_2+\widetilde  F_4^{1/p_3}+\widetilde F_5^\frac{p_5}{p_6}\right)
\cdot 2\left(1+\widetilde F_3^\frac{p_3(2-a)}{p_3(2-a)-1}\right)
=2\mathcal N \Psi_T,
\eeq
where $\mathcal N =1+\widetilde F_1 +\widetilde F_2+\widetilde  F_4^{1/p_3}+\widetilde F_5^\frac{p_5}{p_6}$.

To estimate $\widetilde J$, we calculate the powers of $J_0$ to find that
\beqs
\frac{\alpha_1}{p_1}=\alpha+\lambda(3-2a)-1,\quad 
\frac{\alpha_2}{p_2}=\alpha-\lambda-1,\quad
\alpha=\frac{\alpha_4}{p_3p_4},
\eeqs
and, by the formulas of $p_6$ and $\alpha_5$ in \eqref{ktil} and \eqref{al5}, 
\begin{align*}
\frac{\alpha_5}{p_6}&=\alpha +\frac{p_5(a+\lambda-1)}{p_6(1-a)}=\alpha+\frac{a+\lambda-1}{p_3(2-a)-1}.
\end{align*}
For their minimum  and maximum, we set 
$$h_1'=\max\left\{ -\lambda(3-2a)+1, \lambda+1,0, \frac{-a-\lambda+1}{p_3(2-a)-1} \right\},$$
$$h_2'=\max\left\{\lambda(3-2a)-1,-\lambda-1, 0, \frac{a+\lambda-1}{p_3(2-a)-1}\right\}.$$
Note that
\beqs
 1-\lambda(3-2a)<1<\lambda+1, \quad 
 \frac{1-a-\lambda}{p_3(2-a)-1}<\frac{1-a}{p_3(2-a)-1}<\frac{1-a}{(2-a)-1}=1<\lambda+1.
\eeqs
Thus, $h_1'=\lambda+1=h_1$. Also, it is clear that $h_2'=h_2$.
Hence,
\beqs
0<\alpha-h_1=\alpha-h_1'\le \frac{\alpha_1}{p_1}, 
\frac{\alpha_2}{p_2},
\alpha,
\frac{\alpha_4}{p_3p_4},
\frac{\alpha_5}{p_6}
\le \alpha+h_2'=\alpha+h_2.
\eeqs 
By using inequality \eqref{ee4}, we can estimate $\widetilde J$ as following: 
\beq\label{tildeJfm}
\widetilde J \le 5 (J_0^{\alpha-h_1} +J_0^{\alpha+h_2}).
\eeq

Now, combining  \eqref{I0fm} with \eqref{tildeFfm} and \eqref{tildeJfm}, we obtain
\beqs
\mathcal I_* \le 10C_3\alpha^{\frac{2-a}{1-a}} (1+T) \Big(1+\frac {1}{T_2-T_1}\Big) \Phi_* \mathcal N \Psi_T (J_0^{\alpha-h_1} +J_0^{\alpha+h_2}).
\eeqs
Since $\Phi_* \mathcal N=\mathcal N_1$ and $10C_3=c_{10}$, it follows that
\beq\label{I0i}
\mathcal I_*  \le c_{10}\alpha^{\frac{2-a}{1-a}} (1+T) \Big(1+\frac {1}{T_2-T_1}\Big) \mathcal N_1 \Psi_T (J_0^{\alpha-h_1} +J_0^{\alpha+h_2}).
\eeq
Multiplying \eqref{I0i} by $\alpha/\lambda$  gives
\beqs
B_* \le \frac\alpha\lambda \mathcal I_* 
  \le c_{10} \lambda^{-1}\alpha^{\frac{3-2a}{1-a}}(1+T)\Big(1+\frac {1}{T_2-T_1}\Big) \mathcal N_1 \Psi_T (J_0^{\alpha-h_1} +J_0^{\alpha+h_2})
\eeqs 
which proves \eqref{newS1}.
Dividing \eqref{I0i} by $\widetilde c_0(\alpha-\lambda)/2$ gives
\beqs
A_* \le \frac{2}{\widetilde c_0(\alpha-\lambda)}\mathcal I_* 
 \le 8\mathcal I_*
  \le 8c_{10} \alpha^{\frac{2-a}{1-a}}(1+T)\Big(1+\frac {1}{T_2-T_1}\Big)  \mathcal N_1 \Psi_T (J_0^{\alpha-h_1} +J_0^{\alpha+h_2})
\eeqs
which proves  \eqref{newS2}.
\end{proof}

\begin{lemma}\label{GLk}
Let 
$r_1$ satisfy  \eqref{newro} and denote $r_*$ by \eqref{newrs}.
 Assume the integral in \eqref{wcond1} and
\beq\label{wcond2}
 \int_U  \left(\phi^{-\frac1{1-a}}+\phi^{-\frac{r_1}{1-r_1}}+\phi^{\frac2{r_*}-1}+W_1(x)^{-\frac{r_1}{1-r_1}}\right)\d x
\eeq
are finite. Define
\beq\label{GG}
 \mathcal N_2=\left ( \int_U  \phi(x)^{\frac{2}{r_*}-1}\d x\right )^\frac{r_*}{2}
\left[ \left(\int_U W_1(x)^{-\frac{r_1}{1-r_1}}\d x \right)^{1-r_1}
+\left(\int_U \phi^{-\frac{r_1}{1-r_1}}\d x\right)^{1-r_1}\right]^\frac{1}{r_1}.
\eeq 
Let numbers $\widetilde \kappa>1$ and $p_i>1$ for $i=1,2,\ldots,5$  satisfy \eqref{ktil},  and let $\alpha$ be a number that satisfies \eqref{alp-Large} and 
 \beq\label{Bigalp}
\alpha > \max\left\{\frac{2(\lambda+a-1)}{r_*},\frac{1-\lambda-a}{\widetilde\kappa-1} \right\}.
\eeq
Let $h_1$, $h_2$ be defined by \eqref{h12}, and
\beq\label{kapmax}
    \kappa=\kappa(\alpha)\eqdef 1+\frac {r_*}2+\frac{1-\lambda-a}{\alpha}>1.
\eeq 
If $T>T_2>T_1\ge 0$  then
\beq\label{bfinest2}
\| u\|_{L_\phi^{\kappa\alpha}(U\times (T_2,T))}
\le A_\alpha^\frac 1{\alpha}\Big( \| u\|_{L_\phi^{\widetilde \kappa \alpha}(U\times(T_1,T))}^{\nu_1}+\| u\|_{L_\phi^{\widetilde\kappa \alpha}(U\times(T_1,T))}^{\nu_2}\Big)^\frac 1{\alpha},
\eeq
where  
\beq\label{tilrsdef}
\nu_1=\nu_1(\alpha)\eqdef   \frac{\alpha-h_1}{1+\frac{1}{\alpha(1+r_*/2)}} ,\quad 
\nu_2=\nu_2(\alpha)\eqdef \frac{\alpha+h_3}{1-\frac{\lambda}{\alpha(1+r_*/2)}}
\eeq
with $h_3=\max\{h_2, 1-a-\lambda\}$,
and 
\beq\label{Atildef2}
\begin{split}
A_\alpha &= c_{11} \alpha^{\frac{3(3-2a)}{2(1-a)}}\max\{1,\mathcal N_2\}\left((1+T)\Big(1+\frac1{T_2-T_1}\Big)\mathcal N_1 \Psi_T\right)^{1+r_*/2}
\end{split}
\eeq
with, 
referring to the positive constants $c_7$, $c_8$, $c_9$ in \eqref{ppsi1} and \eqref{c67}, 
\beq  \label{Czero}
c_{11}=2^{4+r_*} \max\{1,c_7^{2-a}\}\max\{c_8,c_9\}^{1+r_*/2}.
\eeq 
\end{lemma}
\begin{proof} 
We apply the parabolic Sobolev inequality  \eqref{ppsi1} in Lemma \ref{WS2} to 
\beqs 
\varphi=\phi, \quad W(x)=W_1(x),\quad 
p = 2 - a,\quad   
s = \lambda+1,
\eeqs
and the interval $[T_2, T]$ in place of $[0,T]$. 
The number $\kappa$ in \eqref{defkappa} takes the specific value in \eqref{kapmax}, while the numbers $\theta_0$ and $m$ in \eqref{tz} and \eqref{mdef} become
\beqs 
\theta_0 =\frac{1}{1+\frac{r_*\alpha}{2(\alpha+1-\lambda-a)}}\in(0,1)
\text{ and }
m=\frac{\alpha+1-\lambda-a}{2-a}\in[1,\alpha).
\eeqs 
The number $\Phi_8$ in \eqref{Jss} turns into $\mathcal N_2$ in \eqref{GG}.
The  conditions in \eqref{alcond} become
\beqs 
\alpha\ge\lambda+1,\ \alpha >\frac{1-a-\lambda}{1-a},\ \alpha>\frac{2(\lambda+a-1)}{r_*}.
\eeqs 
Note that $\frac{1-a-\lambda}{1-a}<1$, then the middle condition can be dropped. The remaining conditions 
 are satisfied thanks to our assumptions \eqref{alp-Large} and  \eqref{Bigalp}.
Then we have from \eqref{ppsi1} and the fact $m<\alpha$ that 
\beq\label{beginest2}
\begin{aligned}
\| u\|_{L_\phi^{\kappa\alpha}(U\times (T_2,T))}
&\le (C_4 \alpha^{\frac 1{r_1}})^\frac{1}{\kappa\alpha} I^{\frac1{\kappa\alpha}}\cdot \essup_{t\in(T_2,T)}\norm{u(\cdot,t)}_{L_\phi^\alpha(U)}^{1-\theta_0},
\end{aligned}
\eeq
where $C_4=c_7^{2-a} \mathcal N_2 $ and 
\beq\label{Iest2}
I=\int_{T_2}^T\int_U  |u(x,t)|^{\alpha+1-\lambda-a}\phi(x)\d x \d t+\int_{T_2}^T\int_U |u(x,t)|^{\alpha-\lambda-1}|\nabla u(x,t)|^{2-a}W_1(x) \d x \d t.  
\eeq
It follows from assumption \eqref{alp-Large}, respectively, \eqref{Bigalp} that 
\beqs 
\alpha>\lambda>a-1+\lambda,\text{ respectively, } 
\alpha>\frac{1-\lambda-a}{\widetilde\kappa-1}.
\eeqs 
These imply $0<\alpha+1-\lambda-a<\widetilde\kappa\alpha$.
Then we can apply H\"older's inequality with the powers $\frac{\widetilde\kappa\alpha}{\alpha+1-\lambda-a}$ and its H\"older's conjugate to the first integral on the right-hand side of \eqref{Iest2}. It results in
\beq\label{Ifirst}
    \int_{T_2}^T\int_U  |u|^{\alpha+1-\lambda-a}\phi\d x \d t
    \le \left(\int_{T_2}^T\int_U |u|^{\widetilde\kappa \alpha}\phi \d x \d t\right)^\frac{\alpha+1-\lambda-a}{\widetilde\kappa\alpha} 
    \left(\int_{T_2}^T\int_U \phi\d x\d t\right)^{1-\frac{\alpha+1-\lambda-a}{\widetilde\kappa \alpha}}.
\eeq

Let $\Phi_*$ be defined by \eqref{Phistar}.
Since $-1/(1-a)<0<1<2/r_*-1$ and by \eqref{ee4}, we in fact have
\beq\label{phiint}
\int_U \phi \d x \le \int_U\left( \phi^{-\frac1{1-a}}+\phi^{\frac2{r_*}-1}\right)\ dx<\infty.
\eeq 
Hence, $\Phi_*\in[1,\infty)$.
With this fact, we  estimate the last integral of \eqref{Ifirst} by
\beq\label{phisimp}
\left(\int_{T_2}^T\int_U \phi\d x\d t\right)^{1-\frac{\alpha+1-\lambda-a}{\widetilde\kappa \alpha}}\le 1+T\int_U\phi \d x\le (1+T)\Phi_*.
\eeq
Denote $\Upsilon= \norm {u}_{L_\phi^{\widetilde \kappa\alpha}(U\times(T_1,T))}$.
We obtain from  \eqref{Iest2},  \eqref{Ifirst} and \eqref{phisimp} and the fact  $\alpha+1-\lambda-a>0$ that
\beq\label{Iest3}
I\le(1+T)\Phi_*\Upsilon^{\alpha+1-\lambda-a}
    +\int_{T_2}^T\int_U |u(x,t)|^{\alpha-\lambda-1}|\nabla u(x,t)|^{2-a}W_1(x) \d x \d t. 
\eeq
To estimate the essential supremum in \eqref{beginest2} and the last double integral in \eqref{Iest3}, we will use Lemma \ref{caccio}.
Note that the requirement $\int_U\phi\d x<\infty$ is satisfied thanks to  \eqref{phiint}.
We compare the powers of $\Upsilon$ in \eqref{newS1}, \eqref{newS2} and \eqref{Iest3}.
Because $a+\lambda-1<\lambda+1= h_1$, it follows that $\alpha-h_1<\alpha-(a+\lambda-1)$. 
Hence,  
$$0<\alpha-h_1\le \alpha+1-\lambda-a,\ \alpha-h_1,\ \alpha+h_2\le \alpha+h_3.$$ 
Applying inequality  \eqref{ee4} yields
\beq \label{UUU}
\Upsilon^{\alpha+1-\lambda-a},\Upsilon^{\alpha-h_1},\Upsilon^{\alpha+h_2}\le \Upsilon^{\alpha-h_1}+ \Upsilon^{\alpha+h_3}.
\eeq 
We denote 
\begin{align}\label{mathcalM}
    \mathcal M_0& =(1+T)\Big(1+\frac1{T_2-T_1}\Big)\mathcal N_1 \Psi_T,\quad
    \mathcal M_1=2\max\{c_8,c_9\}\alpha^{1+\frac{2-a}{1-a}}\mathcal M_0,\\
\label{Sdef}
\mathcal{S}&=\mathcal M_1 (\Upsilon^{\alpha-h_1}+\Upsilon^{\alpha+h_3}),
\end{align}
In calculations below, we very often use the facts
\beq \label{simplef}
\Phi_*,\mathcal N_1,\Psi_T\ge 1\text{  and }c_8,\alpha>1.
\eeq 
In fact, these already imply $\mathcal M_0\ge 1$ and $\mathcal M_1\ge 1$.

It follows from \eqref{newS1}, \eqref{newS2} and \eqref{UUU} that 
\begin{align}\label{Sest1}
&\essup_{t\in(T_2,T)} \int_U u^\alpha(x,t)\phi(x)\d x
\le c_8 \alpha^{\frac{3-2a}{1-a}}\mathcal M_0 (\Upsilon^{\alpha-h_1}+ \Upsilon^{\alpha+h_2})
\le  \mathcal S,\\
\label{Sest2}
&\int_{T_2}^T\int_U |u|^{\alpha-\lambda-1}|\nabla u|^{2-a} W_1(x)\d x \d t 
\le c_9 \alpha^{\frac{2-a}{1-a}}\mathcal M_0 (\Upsilon^{\alpha-h_1}+ \Upsilon^{\alpha+h_2})
\le \alpha^{-1}\mathcal S.
\end{align} 
Also, by \eqref{UUU} and \eqref{Sdef},  
\beq\label{Usim}
\Upsilon^{\alpha+1-\lambda-a}\le \Upsilon^{\alpha-h_1}+ \Upsilon^{\alpha+h_3}=\mathcal M_1^{-1}\mathcal S.
\eeq 
Moreover, one has, thanks to \eqref{mathcalM} and \eqref{simplef}, that
\beq\label{TP}
(T+1)\Phi_*\le (T+1)\mathcal N_1\le \mathcal M_0\le \alpha^{-1-\frac{2-a}{1-a}}\mathcal M_1.
\eeq 
Combining \eqref{Iest3} with \eqref{Sest2}, \eqref{Usim} and  \eqref{TP} yields  
\beq\label{Iest}
I\le ( \alpha^{-1-\frac{2-a}{1-a}}\mathcal M_1)\cdot( \mathcal M_1^{-1}\mathcal S)+ \alpha^{-1}\mathcal S\le \alpha^{-1}\mathcal S+ \alpha^{-1}\mathcal S
= 2\alpha^{-1}\mathcal S.
\eeq
Now, combining \eqref{beginest2},  \eqref{Sest1} and \eqref{Iest}, we obtain
\beq\label{beginest3}
\begin{aligned}
\| u\|_{L_\phi^{\kappa\alpha}(U\times (T_2,T))}
&\le (C_4 \alpha^{\frac 1{r_1}})^\frac{1}{\kappa\alpha} 
(2\alpha^{-1} \mathcal S)^{\frac1{\kappa\alpha}}\cdot \mathcal S^\frac{1-\theta_0}{\alpha}
=\left (2\alpha^{-1+\frac 1{r_1}}C_4 \mathcal S^{\mu_*}\right)^\frac{1}{\kappa\alpha},
\end{aligned}
\eeq
where $\mu_*=1+(1-\theta_0)\kappa$. In fact, one has from \eqref{tk} that 
\beq \label{mustar}
\mu_*=1+r_*/2 .
\eeq 
With the  definition $\mathcal S$ in \eqref{Sdef}, we rewrite \eqref{beginest3} as 
\begin{align*}
\| u\|_{L_\phi^{\kappa\alpha}(U\times (T_2,T))}
&\le \left[ 2\alpha^{-1+\frac 1{r_1}}C_4\mathcal M_1^{\mu_*} (\Upsilon^{\alpha-h_1}+\Upsilon^{\alpha+h_3})^{\mu_*} \right]^\frac{1}{\kappa\alpha} \\
&= \left[ 2\alpha^{-1+\frac 1{r_1}}C_4\mathcal M_1^{\mu_*}\right]^\frac{1}{\kappa\alpha}
 \left[ (\Upsilon^{\alpha-h_1}+\Upsilon^{\alpha+h_3})^\frac{\mu_*}{\kappa} \right]^\frac{1}{\alpha}.
\end{align*}
Applying inequality \eqref{ee1} to estimate $(\Upsilon^{\alpha-h_1}+\Upsilon^{\alpha+h_3})^\frac{\mu_*}{\kappa}$, we obtain
\beq\label{newJ}
\begin{split}
\| u\|_{L_\phi^{\kappa\alpha}(U\times (T_2,T))}
&\le  \left( 2\alpha^{-1+\frac 1{r_1}}C_4 \mathcal M_1^{\mu_*} \right)^\frac{1}{\kappa\alpha}
\left[2^\frac{\mu_*}\kappa \Big(   \Upsilon^{(\alpha-h_1)\frac{\mu_*}\kappa}+\Upsilon^{(\alpha+h_3)\frac{\mu_*}\kappa} \Big)\right]^{\frac 1{\alpha}}\\
&\le  \mathcal M_2^{\frac 1{\kappa\alpha}} \Big( \Upsilon^{\nu_3}+\Upsilon^{\nu_4} \Big)^{\frac 1{\alpha}},
\end{split}
\eeq
where  
     \beq\label{nu12} 
   \nu_3=(\alpha-h_1)\frac{\mu_*}\kappa,\ 
   \nu_4=(\alpha+h_3)\frac{\mu_*}\kappa
    \text{ and  }\mathcal M_2=2^{1+\mu_*}\alpha^{-1+\frac 1{r_1}}\max\{1,C_4\}  \mathcal M_1^{\mu_*}.
   \eeq
Using the facts $\mathcal M_2\ge 1$ and $\kappa>1$, we deduce from \eqref{newJ} a simpler version
\beq\label{newJM}
\| u\|_{L_\phi^{\kappa\alpha}(U\times (T_2,T))}
\le \mathcal M_2^{\frac 1{\alpha}} \Big( \Upsilon^{\nu_3}+\Upsilon^{\nu_4} \Big)^{\frac 1{\alpha}}.
\eeq
By \eqref{kapmax} and \eqref{mustar}, one has
$$\mu_*-\frac\lambda\alpha< \kappa=\mu_*+\frac{1-\lambda-a}{\alpha}< \mu_*+\frac1\alpha,$$
which yields
\beq\label{kmus} 
1+\frac{1}{\alpha \mu_*}> \frac{\kappa}{\mu_*}> 1-\frac{\lambda}{\alpha \mu_*}>0.
\eeq 
(The last inequality comes from the facts $\alpha>\lambda+1>\lambda$ and $\mu_*>1$.)
Using \eqref{kmus} to find lower and upper bounds for $\mu_*/\kappa$ in the formulas of  $\nu_3$ and $\nu_4$ in \eqref{nu12}, and  
noticing also that $\alpha>\lambda+1=h_1$, we obtain 
$$0<\nu_1=\frac{\alpha-h_1}{1+\frac{1}{\alpha(1+r_*/2)}}< \nu_3<\nu_4< \frac{\alpha+h_3}{1-\frac{\lambda}{\alpha(1+r_*/2)}} =\nu_2.$$
Then, by \eqref{ee4},  
\beq\label{U4}
\Upsilon^{\nu_3}\le \Upsilon^{\nu_1} + \Upsilon^{\nu_2}\text{ and }
\Upsilon^{\nu_4} \le \Upsilon^{\nu_1} + \Upsilon^{\nu_2}.
\eeq 
We now rewrite $\mathcal M_2$ more explicitly as
\begin{align*}
    \mathcal M_2
 &     =2^{1+\mu_*}\alpha^{-1+\frac 1{r_1}}\max\{1,c_7^{2-a}\mathcal N_2\}\left (2\max\{c_8,c_9\}\alpha^{1+\frac{2-a}{1-a}}\mathcal M_0\right)^{\mu_*}\\
   & \le  2^{1+2\mu_*}\max\{1,  c_7^{2-a}\} (\max\{c_8,c_9\})^{\mu_*} 
   \alpha^{-1+\frac 1{r_1}+\frac{3-2a}{1-a}\mu_*}\max\{1, \mathcal N_2\}\mathcal M_0^{\mu_*}.
\end{align*}
We calculate and estimate the power of $\alpha$, with the use of the definition \eqref{newrs} of $r_*\in (0,1)$, as follows 
\begin{align*}
&-1+\frac1{r_1}+\left(2+\frac{1}{1-a}\right)\left(1+\frac{r_*}2\right)
=-1+\frac1{r_1}+(2+r_*)+\frac{1}{1-a}\left(1+\frac{r_*}{2}\right)\\
&
=1+\frac1{r_1}+\left (1+\frac{2-a}n-\frac1{r_1}\right )+\frac{1}{1-a}\left(1+\frac{r_*}{2}\right)
=2+\frac{2-a}n+\frac{1}{1-a}\left(1+\frac{r_*}{2}\right)\\
&\le 3+\frac{1}{1-a}\cdot \frac32=\frac{3(3-2a)}{2(1-a)}. 
\end{align*}
For the power of $2$, we quickly have $1+2\mu_*=1+2(1+r_*/2)= 3+r_*$.
Therefore,
\beq\label{lastM1}
\mathcal M_2\le 
2^{3+r_*} \max\{1,  c_7^{2-a}\} (\max\{c_8,c_9\})^{\mu_*} \alpha^\frac{3(3-2a)}{2(1-a)}\max\{1, \mathcal N_2\}  \mathcal M_0^{\mu_*}=\frac{A_\alpha}2.
\eeq 
It follows from \eqref{newJM}, \eqref{U4} and \eqref{lastM1} that
\beqs 
\| u\|_{L_\phi^{\kappa\alpha}(U\times (T_2,T))}
\le \left[ 2\mathcal M_2( \Upsilon^{\nu_1}+\Upsilon^{\nu_2})\right]^{\frac 1{\alpha}}
\le \left[ A_\alpha( \Upsilon^{\nu_1}+\Upsilon^{\nu_2})\right]^{\frac 1{\alpha}}
\eeqs 
which proves \eqref{bfinest2}.
\end{proof}

Next, we apply an extended version of Moser's iteration. We recall the following technical lemma on the upper bounds of certain sequences.

\begin{lemma}[\cite{CHK1}, Lemma A.2]\label{Genn}
Let $y_j\ge 0$, $\kappa_j>0$, $s_j \ge r_j>0$ and $\omega_j\ge 1$  for all $j\ge 0$.
Suppose there is $A\ge 1$ such that
\beqs
y_{j+1}\le A^\frac{\omega_j}{\kappa_j} (y_j^{r_j}+y_j^{s_j})^{\frac 1{\kappa_j}}\quad\forall j\ge 0.
\eeqs
Denote $\beta_j=r_j/\kappa_j$ and $\gamma_j=s_j/\kappa_j$.
Assume
$\bar\alpha \eqdef \sum_{j=0}^{\infty}  \frac{\omega_j}{\kappa_j}<\infty$ 
and the products $\prod_{j=0}^\infty \beta_j$, $\prod_{j=0}^\infty \gamma_j$
 converge to positive numbers  $\bar\beta$, $\bar \gamma$,  respectively.
Then
\beq\label{doub2}
\limsup_{j\to\infty} y_j\le (2A)^{G \bar\alpha} \max\{y_0^{\bar \beta},y_0^{\bar \gamma}\},
\text{ where  }G=\limsup_{j\to\infty} G_j\in(0,\infty)
\eeq
with $G_j=\max\{1, \gamma_{m}\gamma_{m+1}\ldots\gamma_{n}:1\le m\le n< j\}$.
\end{lemma}

In the case $\gamma_j\ge 1$ for all $j\ge 0$, one has
$G_j=\gamma_1\gamma_2\ldots\gamma_{j-1}$ for all $j\ge 2$, hence, the number $G$ in \eqref{doub2} actually is
\beq\label{Gprod}
G=\prod_{j=1}^\infty \gamma_j.
\eeq

\begin{assumption}\label{asmp44}
We assume the following.
\begin{enumerate}[label=\tnum]
    \item Let number $r_1$ satisfy \eqref{newro} and set number $r_*$ by \eqref{newrs}.
    \item Fix a number $\widetilde\kappa\in(1,\sqrt{1+r_*/2})$ and let $p_i>1$, for $i=1,2,\ldots,6$, be as in  \eqref{ktil}.
    \item The integrals in \eqref{wcond1} and \eqref{wcond2} are finite.
\end{enumerate}
\end{assumption}

Recalling  $\mathcal N_1$, $\mathcal N_2$ from \eqref{N1}, \eqref{GG}, respectively, we denote 
\beq \label{N3}
\mathcal N_3=\max\{\mathcal N_1,\mathcal N_2\}.
\eeq 
Then $\mathcal N_3$ involves the weight functions $\phi(x)$, $W_1(x)$, $W_2(x)$ as well as the functions $a_i(x)$ in \eqref{eq2}, while $\Psi_T$ in \eqref{bcquant} involves the boundary data $\psi(x,t)$.

\begin{theorem}\label{LinfU} 
Under Assumption \ref{asmp44}, let $\alpha_0$ be a number such that 
 \beq\label{alp0}
 \begin{split}
\alpha_0>\max\Big\{1+\lambda,\frac{2(\lambda+a-1)}{r_*},\frac{1-\lambda-a}{\widetilde\kappa-1}, \frac{p_1(\lambda(3-2a)-1)}{\widetilde \kappa-p_1}, \\
\frac{p_5(a+\lambda-1)}{(1-a)(\widetilde\kappa -p_6)}, \frac{\lambda+a-1}{1+r_*/2-\widetilde\kappa^2} \Big\}.
\end{split}
\eeq
Then there are positive  constants $\widehat C_0$, $\widetilde  \mu$, $\widetilde \nu$, $\omega_1 $, $\omega_2$, $\omega_3$ independent of $u_0(x)$ and $\psi(x,t)$ such that
if  $T>0$ and $\sigma\in (0,1)$ then
\beq\label{Li1}
\|u\|_{L^{\infty}(U\times(\sigma T,T))} 
\le \widehat C_0 (1+T)^{\omega_1}\Big(1+\frac1{\sigma T}\Big)^{\omega_2}
\mathcal N_3^{\omega_3} \Psi_T^{\omega_2}\max\Big\{ \|u\|^{\widetilde \mu}_{L_\phi^{\widetilde \kappa \alpha_0}(Q_T)},\|u\|^{\widetilde \nu}_{L_\phi^{\widetilde\kappa \alpha_0}(Q_T)}\Big\}.
\eeq
\end{theorem}
\begin{proof}
In order to iterate inequality \eqref{bfinest2} in Lemma \ref{GLk}, we set, for any integer $j\ge 0$,  the numbers $\beta_j=\widetilde \kappa^j\alpha_0$ and $t_j=\sigma T(1-2^{-j})$. 
Thanks to condition \eqref{alp0}, the number $\alpha=\alpha_0$ satisfies \eqref{alp-Large} and \eqref{Bigalp}.
Hence, so do the numbers $\alpha=\beta_j$ for all $j\ge 0$.
We  then can apply estimate \eqref{bfinest2} to $\alpha=\beta_{j}$, $T_2=t_{j+1}$ and $T_1=t_j$ to obtain
\begin{align}
\| u\|_{L_\phi^{\kappa(\beta_{j})\beta_{j}}(U\times (t_{j+1},T))}
&\le   A_{\beta_{j}}^\frac1{\beta_{j}}\Big( \| u\|_{L_\phi^{\beta_{j+1}}(U\times(t_j,T))}^{\widetilde r_j}+\| u\|_{L_\phi^{\beta_{j+1}}(U\times(t_j,T))}^{\widetilde s_j}\Big)^\frac 1{\beta_{j}},\label{ukbb}
\end{align}
where $A_{\beta_j}$ is defined by \eqref{Atildef2}, and, referring to \eqref{tilrsdef}, 
\beq\label{trsj}
\widetilde r_j=\nu_1(\beta_{j})\quad \text{and}\quad\widetilde s_j=\nu_2(\beta_{j}).
\eeq
Since $\widetilde \kappa>1$, the sequence $(\beta_j)_{j=0}^\infty$ is increasing.
From this and \eqref{alp0}, we have 
$$\beta_j\ge \beta_0=\alpha_0>\frac{\lambda+a-1}{1+r_*/2-\widetilde\kappa^2} .$$
Thus, referring to formula \eqref{kapmax} with $\alpha=\beta_j$, one has 
\beq\label{kk}
\kappa(\beta_j)=1+\frac{r_*}{2}+\frac{1-\lambda-a}{\beta_j}
>\widetilde \kappa^2.
\eeq
Note from \eqref{kk} that 
$$\kappa(\beta_{j})\beta_{j}> \widetilde \kappa^2 \beta_{j}=\beta_{j+2}.$$

For $j\ge 0$, define $\mathcal Q_j=U\times (t_j,T)$ and $Y_j=\| u\|_{L_\phi^{\beta_{j+1}}(\mathcal Q_j)}.$
By H\"older's inequality,
\begin{align*}
 Y_{j+1}
 &= \| u\|_{L_\phi^{\beta_{j+2}}(\mathcal Q_{j+1})}
\le\left(T\int_U\phi(x) \ dx\right)^{\frac{1}{\beta_{j+2}}-\frac{1}{\kappa(\beta_{j})\beta_{j}}}\|u\|_{L_\phi^{\kappa(\beta_{j})\beta_{j}}(\mathcal Q_{j+1})}\\
 &\le ((1+T)\Phi_*)^{\frac{1}{\beta_{j+2}}-\frac{1}{\kappa(\beta_{j})\beta_{j}}}\|u\|_{L_\phi^{\kappa(\beta_{j})\beta_{j}}(\mathcal Q_{j+1})}
 \le (1+T)^{\frac{1}{\beta_j}}\Phi_*^{\frac{1}{\beta_{j}}}\|u\|_{L_\phi^{\kappa(\beta_{j})\beta_{j}}(\mathcal Q_{j+1})},
\end{align*}
where $\Phi_*$ is given by \eqref{Phistar}.
Combining this inequality with \eqref{ukbb}, we obtain 
\beq\label{YwithQ}
Y_{j+1}\le \widehat A_j^\frac{1}{\beta_{j}}\big( Y_j^{\widetilde{r}_j}+Y_j^{\widetilde{s}_j}\big)^{\frac 1{\beta_{j}}} \text{ for all }j\ge 0, \text{ where }
\widehat A_j=  (1+T)\Phi_* A_{\beta_{j}}.
\eeq
We estimate $\widehat{A}_j$ now. For the sake of simplicity, we estimate $\Phi_*\le \mathcal N_1$. 
Moreover, from \eqref{N3} and the fact $\mathcal N_1\ge 1$, we have $\max\{1,\mathcal N_2\}\le \mathcal N_3$.
Combining these with  \eqref{Atildef2}  and \eqref{YwithQ} yields 
\begin{align*}
\widehat A_j
&\le  (1+T)\mathcal N_1 \cdot c_{11} \beta_j^{\frac{3(3-2a)}{2(1-a)}}\mathcal N_3\left((1+T)\Big(1+\frac{2^{j+1}}{\sigma T}\Big)\mathcal N_1 \Psi_T\right)^{1+r_*/2} \\
&\le c_{11} 2^{(j+1)(1+r_*/2)}  (1+T)^{2+r_*/2} \Big(1+\frac{1}{\sigma T} \Big)^{1+r_*/2}  \left( \widetilde \kappa^j\alpha_0\right)^{\frac{3(3-2a)}{2(1-a)}} \mathcal N_1^{2+r_*/2}\mathcal N_3
\Psi_T^{1+r_*/2}\\
&\le A_{T,\sigma,\alpha_0}^{j+1},
\end{align*}
where
\beq \label{AAA}
 A_{T,\sigma,\alpha_0}=c_{11} 2^{1+r_*/2}  (1+T)^{2+r_*/2} \Big(1+\frac{1}{\sigma T} \Big)^{1+r_*/2}  \left( \widetilde \kappa\alpha_0\right)^{\frac{3(3-2a)}{2(1-a)}} \mathcal N_1^{2+r_*/2}\mathcal N_3
\Psi_T^{1+r_*/2}.
 \eeq 
Above, we used the properties $\widetilde\kappa>1$, $\alpha_0>1$ and $c_{11}\ge 1$, with the last one coming from \eqref{Czero} and \eqref{simplef}.
Hence, we have 
\beq\label{YY}
Y_{j+1}\le A_{T,\sigma,\alpha_0}^\frac{j+1}{\beta_{j}}\big( Y_j^{\widetilde{r}_j}+Y_j^{\widetilde{s}_j}\big)^{\frac 1{\beta_{j}}} \text{for all }j\ge 0
\eeq

With \eqref{YY}, we apply Lemma \ref{Genn} to the sequence $(Y_j)_{j=0}^\infty$.
We check the convergence of the corresponding sum and products.
Firstly, we have
\beqs
\sum_{j=0}^\infty \frac {j+1} { \beta_j}= \frac 1{ \alpha_0}\sum_{j=0}^\infty \frac {j+1} {\widetilde \kappa^j}=\frac{1}{\alpha_0(1-\widetilde \kappa^{-1})^2}<\infty.
\eeqs
Secondly, by \eqref{trsj},
\beqs
0<\frac{\widetilde r_j}{\beta_j}=\left(1 - \frac{h_1}{\widetilde \kappa^j \alpha_0}\right)\left(1+\frac{1}{\widetilde\kappa^j \alpha_0(1+ r_*/2)}\right)^{-1} <1,  
\eeqs
and
\beqs
\sum_{j=0}^\infty\ln\frac{\widetilde r_j}{\beta_j}=\sum_{j=0}^\infty\ln\left(1 - \frac{h_1}{\widetilde \kappa^j \alpha_0}\right)-\sum_{j=0}^\infty\ln\left (1+\frac{1}{\widetilde\kappa^j \alpha_0(1+ r_*/2)}\right)
\eeqs
which yields $\sum_{j=0}^\infty\ln\frac{\widetilde r_j}{\beta_j}$ is a finite number.
Thirdly, we similarly have from \eqref{trsj} that
\beq\label{q12}
\frac{\widetilde s_j}{\beta_j}=
\left(1+\frac{h_3}{\widetilde\kappa^j\alpha_0}\right)\left(1-\frac{\lambda}{\widetilde\kappa^j\alpha_0(1+r_*/2)}\right)^{-1}>1,
\eeq
and
\beqs 
\sum_{j=0}^\infty\ln{\frac{\widetilde s_j}{\beta_j}}
=\sum_{j=0}^\infty\ln\left(1+\frac{h_3}{\widetilde\kappa^j\alpha_0}\right)
-\sum_{j=0}^\infty\ln \left(1-\frac{\lambda}{\widetilde\kappa^j\alpha_0(1+r_*/2)}\right)
\eeqs 
which  implies $\sum_{j=0}^\infty\ln\frac{\widetilde s_j}{\beta_j}$  is also a finite number.
Therefore, the products 
\beq \label{mnutil}
\widetilde \mu \eqdef   \prod_{j=0}^\infty \frac{\widetilde r_j}{\beta_j}
\text{ and }
\widetilde \nu \eqdef \prod_{j=0}^\infty \frac{\widetilde s_j}{\beta_j}
 \text{ are  positive numbers.}
\eeq 
Then, thanks to \eqref{YY}, \eqref{q12}, \eqref{mnutil} and Lemma \ref{Genn} with the particular case \eqref{Gprod}, we obtain
\beq\label{limY}
\limsup_{j\to\infty} Y_{j}\le (2A_{T,\sigma,\alpha_0})^\omega\max\{Y_0^{\widetilde\mu}, Y_0^{\widetilde\nu}\},
\eeq
where 
\beq \label{oG}
\omega=\mathcal G\sum_{j=0}^{\infty}\frac {j+1}{\beta_j}\text{ with  }
\mathcal G=\prod_{k=1}^\infty (\widetilde s_k/\beta_k)\in(0,\infty).
\eeq 

For $A_{T,\sigma,\varphi}$ in \eqref{AAA},  we simply estimate $\mathcal N_1^{2+r_*/2}\le \mathcal N_3^{2+r_*/2}$.
Therefore, we have
\beq\label{AT2}
(2A_{T,\sigma,\varphi})^\omega
\le C_5\eqdef \widehat C_0 (1+T)^{\omega_1} \Big(1+\frac1{\sigma T}\Big)^{\omega_2}
\mathcal N_3^{\omega_3}  \Psi_T^{\omega_2},
\eeq
where 
$\widehat C_0=\left[2^{2+r_*/2}c_{11} ( \widetilde \kappa \alpha_0)^{\frac{3(3-2a)}{2(1-a)}}\right]^\omega$, and 
\beq \label{omo}
\omega_1=(2+r_*/2)\omega,\quad 
\omega_2=(1+r_*/2)\omega,\quad 
\omega_3=(3+r_*/2)\omega.
\eeq 
Combining \eqref{limY}, \eqref{AT2} with the fact $Y_0=\| u\|_{L_\phi^{\beta_1}(\mathcal Q_0)}=\| u\|_{L_\phi^{\beta_1}(Q_T)}$, we have
\beq\label{limbeta} 
\limsup_{j\to\infty} Y_j
\le C_5  \max\Big\{ \|u\|^{\widetilde \mu}_{L_\phi^{\widetilde \kappa \alpha_0}(Q_T)},\|u\|^{\widetilde \nu}_{L_\phi^{\widetilde\kappa \alpha_0}(Q_T)}\Big\}.
\eeq 
We observe, thanks to H\"older's inequality,  that
\begin{align*}
\|u\|_{L^{\beta_{j}/2}(U\times(\sigma T,T))}&= \left(\int_{\sigma T}^T \int_U \left(|u(x,t)|^{\beta_j/2}\phi(x)^{1/2}\right)\cdot\phi(x)^{-1/2} \d x \d t\right)^{2/\beta_j} \\
&\le \left(\int_{\sigma T}^T \int_U |u(x,t)|^{\beta_j} \phi(x) \d x \d t\right)^{1/\beta_j}
\left (\int_{\sigma T}^T \int_U \phi(x)^{-1} \d x \d t\right)^{1/\beta_j}.
\end{align*}
For $j\ge 1$, simply using the fact $(\sigma T,T)\subset (t_{j-1},T)$  which yields 
$\|u\|_{L_\phi^{\beta_{j}}(U\times(\sigma T,T))} \le Y_{j-1}$, we infer
\beq\label{bbeta2}
\|u\|_{L^{\beta_{j}/2}(U\times(\sigma T,T))}\le Y_{j-1}
\left (T\int_{U} \phi(x)^{-1} \d x\right)^{1/\beta_j}  .
\eeq 
Since $-1/(1-a)<-1<2/r_*-1$, we have from \eqref{ee4} and \eqref{wcond2} that
\beqs 
\int_U \phi^{-1} \d x\le  \int_U  \left(\phi^{-\frac1{1-a}}+\phi^{\frac2{r_*}-1}\right)\d x <\infty.
\eeqs 
By taking the limit superior, as $j\to\infty$, of \eqref{bbeta2} and 
 utilizing \eqref{limbeta},  we obtain 
\begin{align*}
\|u\|_{L^\infty(U\times(\sigma T,T))}
&=\lim_{j\to\infty}\|u\|_{L^{\beta_{j}/2}(U\times(\sigma T,T))} 
\le \limsup_{j\to\infty}Y_{j-1} \\
&\le C_5  \max\Big\{ \|u\|^{\widetilde \mu}_{L_\phi^{\widetilde \kappa \alpha_0}(Q_T)},\|u\|^{\widetilde \nu}_{L_\phi^{\widetilde\kappa \alpha_0}(Q_T)}\Big\}.
\end{align*}
Then the desired estimate \eqref{Li1} follows.
\end{proof}

By combining Theorems \ref{Labound} and \ref{LinfU}, we can establish estimates for  the $L_{x,t}^\infty$-norm of $u(x,t)$ in terms of the initial condition $u_0(x)$  and boundary data $\psi(x,t)$.

\begin{theorem}\label{thm45}
Under Assumption \ref{asmp44}, let $r$ be a number that satisfies \eqref{newrr}  and  $\alpha_0$ be a number that satisfies \eqref{alp0} and 
\beq\label{alzero}
\alpha_0>\frac{2(2-a)(r+a+\lambda-1)}{\widetilde \kappa r_*(1-a)}.
\eeq 
Let $\widetilde  \mu$, $\widetilde \nu$, $\omega_1$, $\omega_2$, $\omega_3$ be the same constants as in Theorem \ref{LinfU}. Denoting $\beta_1=\widetilde \kappa\alpha_0$, let 
\beqs 
V_{0}=1+\int_U \phi(x) u_0(x)^{\beta_1}\d x\text{ and }
M(t) =1+\int_\Gamma (\psi^-(x,t))^{\frac{\beta_1+r}{r}} \d S.
\eeqs
Let the numbers $\mu_{\max}$ and $Z_*$ be as in Theorem \ref{Labound} for $\alpha=\beta_1$.
Assume $T>0$ satisfies 
\beq\label{Tsm0}
\int_0^T M(\tau) \d\tau < \frac{\beta_1}{3 Z_*\mu_{\max} }V_0^{-\frac{\mu_{\max}}{\beta_1}}.
\eeq
Let $\varep$ be any positive number satisfying $\varepsilon<\min\{1,T\}$.
Then one has
\beq\label{Lb1}
\|u\|_{L^{\infty}(U\times(\varepsilon,T))}
\le \widehat C_1\varepsilon^{-\omega_2}(1+T)^{\omega_1} \mathcal N_3^{\omega_3}  \Psi_T^{\omega_2}\max\left\{ \left(\int_0^T\mathcal{V}(t)\d t\right)^{\frac{\widetilde \mu}{\beta_1}}, \left(\int_0^T\mathcal{V}(t)\d t\right)^{\frac{\widetilde \nu}{\beta_1}}\right\}, 
\eeq
where  
\beq\label{newVt}
\mathcal{V}(t)=\left( V_{0}^{-\frac{\mu_{\max}}{\beta_1}}   -\frac{3Z_* \mu_{\max}}{\beta_1} \int_0^t M(\tau) \d\tau  \right)^{-\frac{\beta_1}{\mu_{\max}}}.
\eeq
Consequently, by denoting
\beq\label{Tsm1}
\delta_T\eqdef \frac{3Z_*\mu_{\max}}{\beta_1} V_0^\frac{\mu_{\max}}{\beta_1} \int_0^T M(\tau) \d\tau, 
\eeq
 one has  
\beq\label{Lb2}
\|u\|_{L^{\infty}(U\times(\varepsilon,T))}\le \frac{\widehat C_2 \varepsilon^{-\omega_2}(1+T)^{\omega_1 + \widetilde \nu/\beta_1}}{(1-\delta_T)^{\widetilde \nu/\mu_{\max}}} \mathcal N_3^{\omega_3}  \Psi_T^{\omega_2}   
\Big(1+\norm{u_0}_{L_\phi^{\beta_1}(U)}\Big)^{\widetilde \nu}.
\eeq
Above, $\widehat C_1$ and $\widehat C_2$ are positive constants independent of $T$, $\varep$, $u_0(x)$ and $\psi(x,t)$.
\end{theorem}
\begin{proof}
First of all, applying estimate \eqref{Li1} in Theorem \ref{LinfU} to $\sigma T=\varepsilon<1$, we have 
\beq\label{Li2}
\begin{split}
&\|u\|_{L^{\infty}(U\times(\varepsilon,T))}
\le \widehat C_1 \varepsilon^{-\omega_2}(1+T)^{\omega_1}  \mathcal N_3^{\omega_3}  \Psi_T^{\omega_2}\\
&\quad \times   \max\left\{ \left(\int_0^T\int_U\phi(x)|u(x,t)|^{\beta_1}\d x \d t\right)^{\frac{\widetilde \mu}{\beta_1}},
 \left(\int_0^T\int_U\phi(x) |u(x,t)|^{\beta_1}\d x \d t\right)^{\frac{\widetilde \nu}{\beta_1}}\right\},
 \end{split}
\eeq
where $\widehat C_1=2^{\omega_2}\widehat C_0$.
Next, we use Theorem \ref{Labound} to estimate the last two integrals in \eqref{Li2}.
Thanks to \eqref{alp0}, one has $\beta_1>\alpha_0>\lambda+1$. Together with \eqref{alzero}, we have $\alpha=\beta_1$ satisfies  \eqref{newaa}.
Moreover, thanks to \eqref{Tsm0}, the condition \eqref{tsmall1} is met for $\alpha=\beta_1$. It follows from estimate \eqref{u_wgt_est} of Theorem \ref{Labound} for $\alpha=\beta_1$ that
\beq\label{Lb}
\int_U |u(x,t)|^{\beta_1}\phi(x)\d x \le \mathcal{V}(t) .
\eeq 
Then \eqref{Li2} and \eqref{Lb} imply \eqref{Lb1}.

Thanks to condition \eqref{Tsm0} again, the number $\delta_T$ in \eqref{Tsm1} belongs to the interval $(0,1)$. Note from \eqref{newVt} that $\mathcal V(t)$ is increasing in $[0,T]$. 
Then  one has, for $t\in(0,T)$,
\begin{align*}
\mathcal{V}(t)
&\le \mathcal{V}(T)=\left\{(1-\delta_T)V_0^{-\mu_{\max}/\beta_1}\right\}^{-\beta_1/\mu_{\max}}
=(1-\delta_T)^{-\beta_1/\mu_{\max}} \Big(1+\norm{u_0}_{L_\phi^{\beta_1}(U)}^{\beta_1}\Big)\\
&\le (1-\delta_T)^{-\beta_1/\mu_{\max}}\cdot 2\Big(1+\norm{u_0}_{L_\phi^{\beta_1}(U)}\Big)^{\beta_1}.
\end{align*}
Thus,
\beq \label{Lb3}
\int_0^T \mathcal{V}(t)\d t
\le \frac{ 2(1+T)}{(1-\delta_T)^{\beta_1/\mu_{\max}}}\Big(1+\norm{u_0}_{L_\phi^{\beta_1}(U)}\Big)^{\beta_1}.
\eeq 
Combining \eqref{Lb1} with \eqref{Lb3} and the fact $\widetilde\mu<\widetilde \nu$ yields 
\beqs
\|u\|_{L^{\infty}(U\times(\varepsilon,T))}
\le \widehat C_1\varepsilon^{-\omega_2}(1+T)^{\omega_1} \mathcal N_3^{\omega_3}  \Psi_T^{\omega_2} \left[\frac{ 2(1+T)}{(1-\delta_T)^{\beta_1/\mu_{\max}}}\right]^{\frac{\widetilde \nu}{\beta_1}} 
\Big(1+\norm{u_0}_{L_\phi^{\beta_1}(U)}\Big)^{\widetilde \nu}
\eeqs
which implies the desired estimate \eqref{Lb2} with $\widehat C_2=2^{\widetilde \nu/\beta_1}\widehat C_1$. 
\end{proof}

\begin{example}\label{ex3} 
Same as in Example \ref{ex1}, we consider the case of ideal gas ($\lambda =1/2$) and two-term Forchheimer law ($a=1/2$) for $n=2,3$.  
Under Assumption \ref{asmp44}, the numbers $r_1$ and $r_*$ satisfy \eqref{rr23}, while
\begin{align*}
&\tilde \kappa\in \left(1, \sqrt{15/8-1/(2r_1)} \right) \text{ for } n=2,\text{ or } \tilde \kappa\in \left(1, \sqrt{7/4-1/(2r_1)} \right) \text{ for } n=3.
\end{align*}
Recall that  $1-\lambda-a=0$ and $\lambda(3-2a)-1=0$.
Then the conditions \eqref{alp0} and \eqref{alzero} for $\alpha_0$  become
\beqs
\alpha_0>\max\left\{\frac32,\frac{6r}{\widetilde\kappa r_*}\right\}.
\eeqs
Below, we compute the powers appearing in the estimate \eqref{Lb2}.
In Lemma \ref{caccio}, $h_1=3/2$ and $h_2=0$ which yields $h_3=0$ in  Lemma \ref{GLk}.
From \eqref{q12} and \eqref{mnutil}, the number $\widetilde \nu$   is
\beq\label{specnu}
\widetilde \nu
=\prod_{j=0}^\infty\left(1-\frac{\lambda}{\widetilde\kappa^j\alpha_0(1+r_*/2)}\right)^{-1}
= \begin{cases}\prod_{j=0}^\infty \left(1-\frac{1}{\widetilde\kappa^j\alpha_0(15/4 -1/r_1)}\right)^{-1}& \text {if } n=2,\\
\prod_{j=0}^\infty \left(1-\frac{1}{\widetilde\kappa^j\alpha_0(7/2 -1/r_1)}\right)^{-1} & \text { if } n=3.
\end{cases}
\eeq 
The numbers $\mathcal G$ and $\omega$ in \eqref{oG} are
\beqs
\mathcal G= \begin{cases} \left(1-\frac{1}{\alpha_0(15/4 -1/r_1)}\right) \widetilde \nu & \text {if } n=2,\\
 \left(1-\frac{1}{\alpha_0(7/2 -1/r_1)}\right) \widetilde \nu & \text { if } n=3,
\end{cases}
\text{ and }
\omega=\mathcal G\sum_{j=0}^{\infty}\frac {j+1}{\widetilde\kappa^j\alpha_0} =   \frac{\mathcal G}{\alpha_0(1-\widetilde\kappa^{-1})^2}.
\eeqs
Consequently, the powers $\omega_1$, $\omega_2$, $\omega_3$ given in \eqref{omo} are
\begin{align*}
\omega_1=\frac{ \widetilde \nu}{\alpha_0(1-\widetilde\kappa^{-1})^2}\cdot
\begin{cases} (\frac{23}{8}-\frac{1}{2r_1})\left(1-\frac{1}{\alpha_0(15/4 -1/r_1)}\right) & \text {if } n=2,\\
(\frac{11}{4}-\frac{1}{2r_1})\left(1-\frac{1}{\alpha_0(7/2 -1/r_1)}\right)  & \text { if } n=3,
\end{cases}\\
\omega_2=\frac{ \widetilde \nu}{\alpha_0(1-\widetilde\kappa^{-1})^2}\cdot
\begin{cases} (\frac{15}{8}-\frac{1}{2r_1})\left(1-\frac{1}{\alpha_0(15/4 -1/r_1)}\right)  & \text {if } n=2,\\
(\frac{7}{4}-\frac{1}{2r_1})\left(1-\frac{1}{\alpha_0(7/2 -1/r_1)}\right)  & \text { if } n=3,
\end{cases}\\
\omega_3=\frac{ \widetilde \nu}{\alpha_0(1-\widetilde\kappa^{-1})^2}\cdot
\begin{cases} (\frac{31}{8}-\frac{1}{2r_1})\left(1-\frac{1}{\alpha_0(15/4 -1/r_1)}\right)  & \text {if } n=2,\\
(\frac{15}{4}-\frac{1}{2r_1})\left(1-\frac{1}{\alpha_0(7/2 -1/r_1)}\right)  & \text { if } n=3.
\end{cases}
\end{align*} 
To calculate the number $\mu_{\max}$ appearing in \eqref{Lb2}, we use \eqref{mumax} with $\alpha=\beta_1=\widetilde\kappa\alpha_0$ and the formulas of $r_*$ in \eqref{rr23}. It results in  
\beqs
\mu_{\max}=\frac{3r\beta_1(1+r_*)}{\beta_1r_*-6r}
=3r \beta_1\cdot \begin{cases}  \frac{11/4-1/r_1}{(7/4-1/r_1)\beta_1-6r} & \text {if } n=2,\\
\frac{5/2-1/r_1}{(3/2-1/r_1)\beta_1-6r} & \text { if } n=3.
\end{cases}
\eeqs
In fact, we can have a version of inequality \eqref{Lb2} with even more explicit powers. Indeed, we estimate  $\widetilde\nu$, using the first equation in \eqref{specnu}, by
\begin{align*}
  \widetilde\nu  
  &\le \exp\left\{\sum_{j=0}^\infty \ln \left(1+\frac{\lambda}{\widetilde\kappa^j\alpha_0(1+r_*/2)}\right)\right\}\le \exp\left\{\sum_{j=0}^\infty \frac{\lambda}{\widetilde\kappa^j\alpha_0(1+r_*/2)}\right\}\\
  &= \exp\left\{\frac{1}{\alpha_0(2+r_*)(1-\widetilde\kappa^{-1})}\right\}.
\end{align*}
With this explicit upper bound of $\widetilde\nu$, we can find corresponding upper bounds for $\omega_1$, $\omega_2$, $\omega_3$ and use them in the estimate \eqref{Lb2}. We omit the details. 
\end{example}

\begin{remark}\label{finalrmk}
We remark that one can study problems even more general than \eqref{mainpb}. For example, see Assumption (A1) in \cite{CHK2}, one can  replace condition \eqref{Zb} for $Z(u)$ with
\beqs
|Z(u)|\le d_0 |u|^{\ell_Z}.
\eeqs
A source/sink term $f(x,t,u)$ can be added to the right-hand side of the first equation in \eqref{mainpb} with 
\beqs
|f(x,t,u)|\le f_1(x,t)+f_2(x,t)|u|^{\ell_f}.
\eeqs
Moreover, the boundary term $\psi (x,t)u^\lambda$ in \eqref{mainpb} can be replaced with $-B(x,t,u)$, where $B(x,t,u)$ satisfies
\beqs
B(x,t,u)\le \varphi_1(x,t)+\varphi_2(x,t) |u|^{\ell_B}.
\eeqs
Above, constants $d_0$, $\ell_Z$ are positive, constants $\ell_f$, $\ell_B$ are non-negative and functions $f_1$, $f_2$, $\varphi_1$, $\varphi_2$ are non-negative.
Then we can combine the methods in this paper and \cite{CHK2} to obtain similar results to Theorems \ref{Labound} and \ref{thm45}. However, the calculations will be even more complicated.
\end{remark}

\medskip
\noindent\textbf{Data availability.} 
No new data were created or analyzed in this study.

\medskip
\noindent\textbf{Funding.} No funds were received for conducting this study. 

\medskip
\noindent\textbf{Conflict of interest.}
There are no conflicts of interests.

\myclearpage

\bibliography{paperbaseall}{}
\bibliographystyle{abbrv}
\end{document}